\documentclass[11pt,openbib]{article}
\usepackage{amsmath}
\usepackage{amssymb}
\usepackage{amsfonts}
\usepackage[T1]{fontenc}
\usepackage{graphicx}

\newcommand{\R}{\mathbb{R}}
\newcommand{\C}{\mathbb{C}}

\newcommand{\Z}{\mathbb{Z}}
\newcommand{\dee}{\mathrm{d}}
\newcommand{\comp}{\, \raisebox{2pt}{$\scriptstyle\circ \, $}}
\newcommand{\rowspace}{\rule{0pt}{16pt}}
\newcommand{\onehalf}{\mbox{$\frac{\scriptstyle 1}{\scriptstyle 2}\,$}}
\newcommand{\ttfrac}[2]{\mbox{${\scriptstyle \frac{{#1}}{{#2}}}$}}
\newcommand{\spann}{\mathop{\rm span}\nolimits}
\newcommand{\lefthook}{\mbox{$\, \rule{8pt}{.5pt}\rule{.5pt}{6pt}\, \, $}}

\newcommand{\smalldbydt}{\mbox{${\scriptstyle \frac{\dee}{\dee t}}
\rule[-6pt]{.5pt}{12pt} \raisebox{-6pt}{$ \, {\scriptscriptstyle t=0}$}$}}
\newcommand{\setrule}{\, \rule[-4pt]{.5pt}{13pt}\, }
\newcommand{\bigspace}{\bigskip \par\noindent}
\newcommand{\ad }{\mathop{\mathrm{ad} }\nolimits}

\newenvironment{proof}[1][Proof]{\noindent\textbf{#1.} }{\ \rule{0.5em}{0.5em}}
\allowdisplaybreaks

\begin{document}

\title{Bohr-Sommerfeld-Heisenberg Quantization \\
of the Mathematical Pendulum\footnotemark}
\author{Richard Cushman and J\k{e}drzej \'{S}niatycki\\ 
Department of Mathematics and Statistics\\
University of Calgary, Calgary, Alberta, Canada.
\date{}}
\maketitle
\footnotetext{Appeared in Journal of Geometric Mechanics \textbf{10} (2018) 419--443.}
\footnotetext{Corrected: 21 November 2019}
\begin{abstract}
In this paper we give the Bohr-Sommerfeld-Heisenberg quantization of the 
mathematical pendulum.
\end{abstract}

\section{Introduction}
%%%%%%%%%%%%%%%

The Dirac's formulation of quantum mechanics \cite{dirac25} can be described 
as a precursor of the theory of ${\C }^{\ast }$ algebras. Quantum 
observables are self adjoint operators in a complex vector space of quantum states. In 
chapter 3 of \cite{dirac30} Dirac represents quantum states as functions on the spectrum 
of the maximal abelian subalgebra (complete set of commuting observables). Classical 
Hamiltonian mechanics may be regarded as the limit of quantum mechanics when $\hbar $ tends to zero. There are quantum systems without classical analogues. 
\par Quantization is an attempt to find a quantum system corresponding to a given 
classical system. Since there may be several different approaches, quantization may give 
inequivalent results. Because quantum observables may be represented as operators on 
the space of functions on the spectrum of the maximal abelian subalgebra, the usual 
approach to quantization is to identify a complete set of commuting observables and to 
study operators on its spectrum.
\par For a completely integrable Hamiltonian system, 
Bohr-Sommerfeld quantization \cite{bohr,sommerfeld} of the action variables gives rise to a space of quantum states and a complete set of commuting observables acting of this space of states. Bohr-Sommerfeld theory does not provide operators of transition between the 
eigenstates of operators corresponding to the actions. These transitions are accounted 
for by shifting operators. Because the general theory of these operators requires an extension of geometric quantization to locally Hamiltonian vector fields, which is far a field from the topic of this paper, we refer the reader to \cite{cushman-sniatycki18}. However, we do treat a special case relevant to this paper in the appendix. The commutation relations satisfied by the shifting operators are the same as the commutation relations satisfied by formal quantization of the functions 
${\mathrm{e}}^{\pm i \vartheta }$, where $\vartheta $ is an angle in the action angle coordinates for the integrable system. Moreover, if $\vartheta $ were a single-valued function, then its Hamiltonian vector field $X_{\vartheta }$ would generate a local group 
${\mathrm{e}}^{tX_{\vartheta }}$ of local symplectomorphisms of the phase space
preserving the Bohr-Sommerfeld polarization, which would lift to a local
group ${\mathrm{e}}^{tZ_{\vartheta }}$ of local quantomorphisms. Since the
angle $\vartheta $ is a multi-valued function, ${\mathrm{e}}^{tZ_{\vartheta }}$
is not well defined for $t\neq nh$, where $h$ is Planck's constant and $n\in 
\mathbb{Z}.$ However, the shifting operators, given by 
${\mathrm{e}}^{\pm hZ_{\vartheta }}$ are well defined and correspond to the operators of
multiplication by $\mathrm{e}^{\pm i\vartheta }$. The existence of shifting operators answers Heisenberg's criticism \cite{heisenberg25} of the Bohr-Sommerfeld theory. 
\par In geometric quantization, a complete set of commuting observables corresponds to 
a polarization. For a completely integrable Hamiltonian system with a regular foliation by 
Lagrangian tori, we get Bohr-Sommerfeld theory by choosing a polarization tangent to the 
tori of the foliation \cite{cushman-sniatycki12}. Taking into account the existence of shifting operators, we obtain a \emph{full} geometrically based quantum theory. We do not try to compare the results of our quantization scheme with observations. For readers who would like to compare the energy spectra of the Schrodinger and the Bohr-Sommerfeld quantizations of the mathematical pendulum, we provide implicit equations for the energy spectrum in Bohr-Sommerfeld theory. In interesting completely integrable systems \cite{cushman-sniatycki13}, the foliation by tori is not regular and we have to take into account the singularities of the polarization to obtain the Bohr-Sommerfeld quantum spectrum. 
\par In this paper we discuss how to treat the singularities in the mathematical pendulum. 

\section{The classical mathematical pendulum}\label{section-2}
%%%%%%%%%%%%%%%%%%%%% 

\subsection{The basic setup}\label{subsection-2.1}
%%%%%%%%%%%%%%%%%

We consider the classical mathematical pendulum, which is a Hamiltonian
system on $T^{\ast }S^{1}=\mathbb{R}\times S^{1}=\mathbb{R}\times \left( 
\mathbb{R}/2\pi \mathbb{Z}\right) $, the cotangent bundle of the circle $%
S^{1}$, with coordinates $(p,\alpha )$, symplectic form $\omega =\mathrm{d}
p\wedge \mathrm{d} \alpha $, and $1$-form $\theta = p \, \dee \alpha $. The Hamiltonian of the system is 
\begin{equation}
H:T^{\ast }S^{1}\rightarrow \mathbb{R}:(p,\alpha )\mapsto %
\mbox{$\frac{\scriptstyle 1}{\scriptstyle 2}\,$} p^{2}-\cos \alpha +1.
\label{eq-s1one}
\end{equation}
The Hamiltonian vector field $X_{H}$ of $H$ satisfies $X_{H} 
\mbox{$\,
\rule{8pt}{.5pt}\rule{.5pt}{6pt}\, \, $} (\mathrm{d} p\wedge \mathrm{d}
\alpha )=-p\, \mathrm{d} p-\sin \alpha \, \mathrm{d} \alpha $ so that 
\begin{equation}
X_{H}(p,\alpha )=-\sin \alpha \,\frac{\partial }{\partial p}+p\, \frac{%
\partial }{\partial \alpha }.  \label{eq-s1two}
\end{equation}
Its integral curves are solutions of Hamilton's equations 
\begin{align}
\frac{\mathrm{d} p}{\mathrm{d} t}& =-\sin \alpha \quad \text{and} \quad 
\frac{\mathrm{d} \alpha }{\mathrm{d} t} =p.  \label{eq-s1twostar}
\end{align}
The Hamiltonian $H$ has two critical points: one at $(0,0)$ with $H(0,0)=0$
and the other at $(0,\pi )$ with $H(0,\pi )=2$. These correspond to a stable
elliptic and an unstable hyperbolic equilibrium point of $X_{H}$,
respectively.
\mbox{}
\vspace{.1in}\begin{tabular}{l} 
\vspace{-.4in} \\
\setlength{\unitlength}{2pt} \\
\includegraphics[width=350pt]{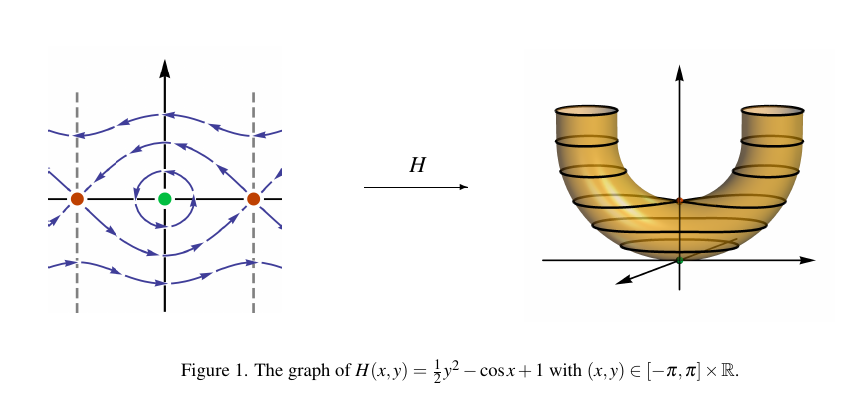}\\
\vspace{-.45in}
\end{tabular}

\subsection{Action-angle coordinates}\label{subsection-2.2}
%%%%%%%%%%%%%%%%%%%

In this subsection we find action-angle coordinates $(I, \vartheta )$ for the mathematical
pendulum. \medskip 

First, we introduce the action function $I$ on $T^{\ast }S^{1}$ such
that for every connected component $C(e)$ of the energy level $H^{-1}(e)$,\ the restriction of 
$I$ to $C(e)$ is 
\begin{equation}
I(e) = I|_{C(e)}=  \ttfrac{1}{2\pi } \int_{C(e)} \theta 
= \ttfrac{1}{2\pi } \int_{C(e)}p\,\mathrm{d}\alpha .  
\label{0}
\end{equation}
Before giving explicit expressions for $I$ and $\vartheta $ we compute the Poisson bracket 
$\{ I, \vartheta \} $ as follows:
\begin{align}
\{ I, \vartheta \} & = L_{X_{\vartheta }}I = \ttfrac{1}{2\pi} \int_{C(e)} L_{X_{\vartheta }} \theta 
= \ttfrac{1}{2\pi}  \int_{C(e)} [ X_{\vartheta } \, \lefthook \, \dee \theta + \dee  (X_{\vartheta } \, 
\lefthook \, \theta ) ] \notag \\
& = - \ttfrac{1}{2\pi}  \int_{\gamma } \dee \vartheta  = -1, \label{eq-1}
 \end{align}
since $\omega = \dee \theta $ and $C(e)$ is parametrized by a periodic integral curve $\gamma $ of $X_H$ of period $T = T(e)$. We reparametrize $C(e)$ using $\vartheta = \frac{2\pi}{T}\, t$, which is 
the angle function. Because the matrix of the symplectic form $\omega $  in action angle coordinates is 
\begin{displaymath}
\mbox{\small $\left( \begin{pmatrix} 0 & \{ I, \vartheta \} \\ \{ \vartheta , I \} & 0 \end{pmatrix}^{\! \! -1} \right)^t = 
\begin{pmatrix} 0 & -1 \\ 1 & 0 \end{pmatrix} $}, 
\end{displaymath}
it follows that $\omega = \dee I \wedge \dee \vartheta $. Similarly, the Poisson bracket 
$\{ I, H \} $ is computed as follows: 
\begin{align}
\{ I, H \} & = L_{X_H} I = \ttfrac{1}{2\pi} \int_{C(e)} L_{X_H} \theta = 
 \ttfrac{1}{2\pi}  \int_{C(e)} [ X_H \, \lefthook \, \dee \theta + \dee  (X_H  \, 
\lefthook \, \theta ) ] \notag \\
& = \ttfrac{1}{2\pi} \int_{\gamma } \dee ( -H + X_H \, \lefthook \, \theta ) = 0, \notag 
\end{align}
since the curve $\gamma $ is closed. Thus $I$ is constant on the integral curves of $X_H$. So 
$I$ is constant on $C(e)$. Consequently, 
\begin{equation}
\ttfrac{1}{2\pi} \int_{C(e)} I \, \dee \vartheta = \ttfrac{1}{2\pi} I|_{C(e) } 
\int_{\gamma } \dee \vartheta = I(e). 
\label{eq-s2onenew}
\end{equation}

We now give explicit expressions for the action $I$ and the angle $\vartheta $ of the mathematical 
pendulum. There are two cases. \medskip 

\noindent \textsc{Case} 1. $0 < e < 2$.

\noindent We denote by $I_{0}$ the restriction of $I$ to the region $P_{0}= \{(p,\alpha )\in T^{\ast }S^{1}\mid H(p,\alpha )<2\} $. Because $(0,0)$ is a nondegenerate minumum of the Hamiltonian $H$ with
minimum value $0$, for $e$ near $0$ the level set $H^{-1}(e)$ is
diffeomorphic to a circle $S^{1}$ and hence is connected. From the Morse isotopy lemma it follows
that for every $e$ with $0<e<2$ the level set $H^{-1}(e)$ is diffeomorphic
to a circle and hence is connected. By definition 
\begin{equation}
I_{0}(e)= \ttfrac{1}{2\pi} \int_{H^{-1}(e)}p\,\mathrm{d}\alpha =
\ttfrac{1}{\pi } \int_{{\alpha }^{-}}^{{\alpha }^{+}}\sqrt{2\left( e-(1-\cos \alpha )\right) }\,\mathrm{d}\alpha ,
\label{eq-s1three}
\end{equation}
where $e=1-\cos {\alpha }^{\pm }$, which implies that ${\alpha }^{-}=-{%
\alpha }^{+}$, since $\cos $ is an even function. Therefore 
\begin{equation}
I_{0}(e)=\ttfrac{4}{\pi} e\int_{0}^{\pi /2}\frac{{\cos }^{2}\varphi }{\sqrt{1-\frac{e}{2}{%
\sin }^{2}\varphi }}\,\mathrm{d}\varphi ,
\label{eq-s1threedot}
\end{equation}
using the identity $\cos \alpha = 1 -2 \, {\sin }^2\frac{\alpha}{2} $ and the change of variables 
$\sin \frac{\alpha }{2} = \sin \frac{{\alpha }^{+}}{2} \sin \varphi $. We check some
limiting cases. First when $e\nearrow 2$ we obtain 
$\lim_{e\nearrow2}I_{0}(e)=\frac{8}{\pi }$. When $e\searrow 0$ we find that 
$I_{0}(e) \sim \frac{4e}{\pi} \,\int_{0}^{\pi /2}{\cos }^{2}\varphi \,\mathrm{d}\varphi = e$, 
which is what is given by
the harmonic oscillator. \medskip

\noindent We now find the corresponding angle $\vartheta _{0}$. By
definition 
\begin{equation*}
\vartheta _{0}=\frac{{\small 2\pi }}{{\small T}}t=\frac{{\small 2\pi }}{{\small T}}\int_{-{\alpha }^{+}}^{\alpha }\frac{\mathrm{d}\alpha }{\sqrt{2(e-(1-\cos \alpha ))}}=\frac{{\small 4\pi }}{{\small T}}\,\int_{0}^{\varphi
}\frac{1}{\sqrt{1-\frac{e}{2}{\sin }^{2}\varphi }}\,\mathrm{d}\varphi ,
\end{equation*}
where $T=T(e)$ is the period of the motion of the mathematical pendulum on $%
H^{-1}(e)$. From Hamilton's equations it follows that 
\begin{align}
T& = 4\int_{0}^{\pi /2}\frac{1}{\sqrt{1-\frac{e}{2}{\sin }^{2}\varphi }}\, \mathrm{d}\varphi .  \label{4}
\end{align}
Again we check some limiting cases. First, when $e\nearrow 2$ we
find that $T\nearrow \infty $. So $\vartheta _{0}\searrow 0$. Second, when $%
e\searrow 0$ we get $T\searrow 4\int_{0}^{\pi /2} \dee \varphi =2\pi $. So $%
\vartheta _{0}\searrow 2\varphi =\alpha $, which checks with the angle given
by the harmonic oscillator. \medskip

\noindent \textsc{Case} 2. $e > 2$.%

\noindent First we find the restrictions $I_{\pm }$ of $I$ to the regions 
$ P_{\pm }=\{(p,\alpha )\in T^{\ast }S^{1}\mid H(p,\alpha )>2\text{, \ }\pm
p>0\} $. Because $(0,\pi )$ is a nondegenerate critical point of Morse index $1$ of
the Hamiltonian $H$ with critical value $2$, for $e>2$ but near to $2$ the
level set $H^{-1}(e)$ is diffeomorphic to the disjoint union of two circles $%
C_{\pm }(e)$. By the Morse isotopy lemma it follows that for all $e>2$ the
level set $H^{-1}(e)$ is diffeomorphic to $C_{-}(e)\coprod C_{+}(e)$. By
definition 
\begin{align}
I_{\pm }(e)& = \ttfrac{1}{2\pi} \int_{C_{\pm }(e)}p\,\mathrm{d}\alpha = 
\ttfrac{1}{\pi} \int_{-\pi }^{\pi }\sqrt{2(e-(1-\cos \alpha ))}\,\mathrm{d}\alpha  \notag \\
& = \ttfrac{4\sqrt{2e}}{\pi } 
\int_{0}^{\pi /2}\sqrt{1-\ttfrac{2}{e}\,{\sin }^{2}\varphi }\,\,\mathrm{d}\varphi . 
\label{eq-s1threedagg}
\end{align}%
We check two limiting cases. When $e\searrow 2$, $\lim_{e\searrow 2}I_{\pm }(e)=
\ttfrac{4}{\pi} \int_{0}^{\pi /2}\cos \varphi \,\mathrm{d}\varphi =\ttfrac{4}{\pi} $, 
which is one half of the action $I(e)$ at $e=2$. This is
correct because as $e\searrow 2$ the component $C_{\pm }(e)$ of $H^{-1}(e)$
converges to $H^{-1}(2)\cap \{\pm p\geq 0\}$. When $e\nearrow \infty $, we
get $I_{\pm }(e) \sim \sqrt{2e}$. \medskip

\noindent We now find the corresponding angle ${\vartheta }_{\pm}$. By definition 
\begin{equation*}
{\vartheta }_{\pm }=\frac{2\pi }{T_{\pm}}t=
\frac{2\pi }{T_{\pm }}\int_{-\pi
}^{\alpha }\frac{\mathrm{d}\alpha }{\sqrt{2(e-(1-\cos \alpha ))}}=\frac{2\pi 
}{T_{\pm }}\sqrt{\ttfrac{2}{e}}\,\int_{0}^{\varphi }\frac{1}{\sqrt{1-\frac{2}{e%
}{\sin }^{2}\varphi }}\,\mathrm{d}\varphi .
\end{equation*}%
where $T_{\pm }=T_{\pm }(e)$ is the period of the motion of the mathematical
pendulum on $H^{-1}(e)$. From Hamilton's equations it follows that 
\begin{align}
T_{\pm }& =\int_{-\pi }^{\pi }\frac{\mathrm{d}\alpha }{\sqrt{2(e-(1-\cos \alpha ))}}   
=\sqrt{\frac{2}{e}}\,\int_{0}^{\pi/2}\frac{1}{\sqrt{1-\frac{2}{e}\,{\sin }^{2}\varphi }}\,\mathrm{d}\varphi .
\label{7}
\end{align}
Again we check some limiting cases. First, when $e\searrow 2$ we
find that $T_{\pm }\nearrow \infty $. So ${\vartheta }_{\pm }\searrow 0$.
Second, when $e\nearrow \infty $ we get $T_{\pm }\sim \frac{\pi }{\sqrt{2e}}$. 
So ${\vartheta }\sim 4\varphi =2\alpha $. \medskip

\noindent It follows from the above discussion that the action function $I$, defined
by equation (\ref{0}) is continuous on $[0, \infty )$. However, $I(e)$ is not smooth at $e=2$, 
see Dullin \cite{dullin}. 

\section{Elements of geometric quantization}\label{section-3}
%%%%%%%%%%%%%%%%%%%%%%%

In this section, we review the elements of geometric quantization applicable
to the mathematical pendulum following \cite{sniatycki80}.\medskip

\noindent Consider a trivial complex line bundle $L=\mathbb{C}\times T^{\ast }S^{1}$
with projection map $\rho :L\rightarrow T^{\ast }S^{1}: \big( z,(p,\alpha) %
\big)\mapsto (p,\alpha )$ and trivializing section $\lambda
_{0}:T^{\ast}S^{1}\rightarrow L:(p,\alpha ) \mapsto (1,(p,\alpha ))$. Define
a connection $\nabla $ on $L$ by setting 
\begin{equation}
\nabla \lambda _{0}=-i\hbar ^{-1}\theta \otimes \lambda _{0},  
\label{nabla}
\end{equation}
where $\hbar $ is Planck's constant divided by $2\pi $ and $\theta =p\, 
\mathrm{d} \alpha $ is the canonical $1$-form on $T^{\ast }S^{1}$. Since 
$\omega =\mathrm{d} \theta $, it follows that the curvature of the connection 
$\nabla $ is $\frac{-1}{2\pi \hbar }\, \omega $. \medskip

\noindent We \ consider the geometric quantization of the mathematical pendulum with
respect to the singular polarization $D$ of $T^{\ast }S^{1}$
consisting of all integral curves of the Hamiltonan vector field 
$X_{H}$ (\ref{eq-s1two})\footnote{Throughout this paper we will use the shorthand 
$\mathit{mathsymbol} \, (\mathit{Number})$ to mean $\mathit{mathsymbol}$ given in 
equation $(\mathit{Number})$. For example, $X_H \, (2)$ means $X_H$ given in 
equation $(2)$.} associated to the Hamiltonian function $H$ (\ref{eq-s1one}).
This means that quantum states of the mathematical pendulum are represented by sections of
the prequantization line bundle $L$ that are covariantly constant along $D$. It should be noted that this representation is not unique. Multiplication of every section by a constant phase factor leads to an equivalent representation. For $e\notin \{0,2\},$ the leaves of $D$ are smooth and are topological circles due to the conclusions of subsection 2.2. Moreover, this polarization has singularities consisting of the equilibrium points $(0,0)$ and $(0,\pi )$ of $X_{H}$ and two
homoclinic orbits of $X_{H}$, which have $(0,\pi )$ as a common boundary.
Therefore, we are extending geometric quantization to a \emph{singular} 
polarization, which leads to the difficulties encountered here. \medskip

\section{Bohr-Sommerfeld conditions}\label{section-4}
%%%%%%%%%%%%%%%%%%%%%%

Consider an integral curve $\gamma :\mathbb{R}\rightarrow T^{\ast}S^{1}:t\mapsto 
\gamma (t)=\big(p(t),\alpha (t)\big)$ of the Hamiltonian
vector field $X_{H}$. Suppose that $e=H(\gamma (t))$ is not $0$ or $2.$
Then, $\gamma $ is periodic with period $T\neq 0$. The cases when $e=0$ and $%
e=2$ will be discussed separately. \medskip

\noindent Let $\sigma :T^{\ast }S^{1}\rightarrow L$ be a section of the
prequantization line bundle that is covariantly constant along $D$. Then ${\gamma}^{\ast }\sigma :$ $\mathbb{R}\rightarrow L$ is a horizontal lift of $\gamma $ to $L$. ${\gamma}^{\ast }\sigma $ is periodic with period $T$ if either the restriction of the connection $\nabla $ to the image of $\gamma $ has trivial holonomy group, or $\sigma $ restricted to the image of $\gamma $
is identically zero. \medskip 

\noindent \textbf{Theorem 4.1}
\label{BS}\textit{Let $\gamma :[0,T]\rightarrow T^{\ast }S^{1}$ be a periodic
integral curve of $X_{H}$ with period $T\neq 0$. The holonomy group of $\nabla $,
restricted to the image $\mathrm{im}\, \gamma $ of $\gamma $, is trivial if and only if the action
integral} 
\begin{equation}
I_{\gamma }=I|_{\mathrm{im}\, \gamma } = 
\ttfrac{1}{2\pi} \int_{0}^{T}\gamma ^{\ast }\theta  \, \, \dee t =n\hbar , 
\label{Bohr-Sommerfeld}
\end{equation}%
\textit{for some $n \in \Z $.}

\begin{proof}
Consider an integral curve $\gamma :\mathbb{R}\rightarrow T^{\ast
}S^{1}:t\mapsto \gamma (t)=\big(p(t),\alpha (t)\big)$ of the Hamiltonian
vector field $X_{H}$. It satisfies Hamilton's equations (\ref{eq-s1twostar}%
). Suppose that $e=H(\gamma (t))\in (0,2)$. Then, the curve $\gamma $ is
periodic with period $T$, see (\ref{4}). Let 
\begin{equation*}
\widetilde{\gamma }:[0,T]\rightarrow L:t\mapsto \big(z(\gamma (t)),\gamma (t)%
\big)=z( \gamma (t)) \,{\lambda }_{0}(\gamma (t)) 
\end{equation*}%
be a horizontal lift of $\gamma $. Then the covariant derivative $\frac{D}{%
\mathrm{d}t}\widetilde{\gamma }(t)$ of $\widetilde{\gamma }$ must vanish.
Equation (\ref{nabla}) implies that 
\begin{eqnarray*}
\frac{D}{\mathrm{d}t}\widetilde{\gamma }(t) &=&\frac{D}{\mathrm{d}t}\big(
z(\gamma (t))\lambda _{0}(\gamma (t))\big) \\
&=&\frac{\mathrm{d}z}{\mathrm{d}t}\,\lambda _{0}(\gamma (t))-i\hbar
^{-1}z\langle \theta \mid X_{H}\rangle (\gamma (t))\,\lambda
_{0}(\gamma (t)) \\
&=&\Big( \frac{\mathrm{d}z}{\mathrm{d}t}-i\hbar
^{-1}z\big\langle p\,\mathrm{d}\alpha \mid -\sin \alpha \,\frac{\partial }{\partial p}
+p\,\frac{\partial }{\partial \alpha }\big\rangle \Big) {\lambda}_{0}(\gamma (t)) \\
&=&\Big( \frac{\mathrm{d}z}{\mathrm{d}t}-i\hbar ^{-1}p(t)^{2}\,z \Big)
\lambda _{0}(\gamma (t)),
\end{eqnarray*}
where $\mbox{$\frac{\scriptstyle 1}{\scriptstyle 2}\,$}p(t)^{2}-\cos \alpha
(t)+1=e$. Hence, $p(t)=\pm \sqrt{2(e-(1-\cos \alpha (t)))}$. \linebreak Because 
$\frac{\mathrm{d}z}{\mathrm{d}t}=\frac{\mathrm{d}z}{\mathrm{d}\alpha }\frac{%
\mathrm{d}\alpha }{\mathrm{d}t}=\frac{\mathrm{d}z}{\mathrm{d}\alpha }\,p$,
the curve $\widetilde{\gamma }$ is horizontal (= covariantly constant) 
\rule{0pt}{12pt}if $%
\frac{\mathrm{d}z}{\mathrm{d}\alpha }\,p(t)-i\hbar ^{-1}p(t)^{2}z= p(t)\big( \frac{\mathrm{d}z}{\mathrm{d}\alpha } 
-i\hbar ^{-1}p(t)z \big) = 0$, that is, 
\begin{equation}
-i\hbar \,\frac{1}{z}\frac{\mathrm{d}z}{\mathrm{d}\alpha }=\pm \sqrt{%
2(e-(1-\cos \alpha ))}.  \label{paralleltransport}
\end{equation}%
Here the $+$ sign corresponds to $\alpha \in \lbrack 0,\alpha ^{+}]$ and the 
$-$ sign corresponds to $\alpha \in \lbrack \alpha ^{-},0]=[-\alpha ^{+},0]$.
Integrating (\ref{paralleltransport}) from $\alpha ^{-}$ to $\alpha ^{+}$
and using the fact that $\cos $ is an even function, we get 
\begin{equation*}
-i\hbar \ln \left( \frac{z(\alpha ^{+})}{z(\alpha ^{-})}\right)
=2\int_{\alpha ^{-}}^{\alpha ^{+}}\sqrt{2(e-(1-\cos \alpha ))}\,\mathrm{d}\alpha =2\pi I(e)
\end{equation*}%
by equation (\ref{eq-s2onenew}). The horizontal lift $\widetilde{\gamma }$ of the closed curve 
$\gamma $ is a closed curve in the line bundle $L$ if and only if 
$z({\alpha }^{+})=z(\alpha ^{-})$. Since $\ln $ is a multivalued function and $\ln
1=2\pi ni$, it follows that $\widetilde{\gamma }$ is a closed curve in $L$
if and only if we have $I|_{C(e)}=n{\hbar} $. \medskip  

\noindent For $e>2$, we have $H^{-1}(e)=C_{-}(e)\amalg C_{+}(e)$, and there are integral curves $%
\gamma ^{-}$ and $\gamma ^{+}$ of $X_{H}$ such that where $C_{-}(e)$ is the 
image of $\gamma ^{-}$ and $C_{+}(e)$ are image of $\gamma ^{+}$. The same
argument as in the preceding paragraph shows that the horizontal lift $\tilde{\gamma}^{-}$ of $%
\gamma ^{-}$ is a closed curve in $L$ if and only if $I_{-}(e)=m_{-}\hbar $, 
where $m_{-}$ is an integer. Similarly, the horizontal lift ${\widetilde{\gamma}}^{+}$ of 
$\gamma ^{+}~$is a closed curve in $L$ if and only if $I_{+}(e)=m_{+}\hbar ,$ where $m_{+}$ is an integer. Since $I_{-}(e)=I_{+}(e)$, it follows that $m_{-}=m_{+}$. Moreover, 
$I_{-}(e)=I|_{C_{-}(e)}$ and $I_{+}(e)=I|_{C_{+}(e)}$. 
\end{proof}
\bigspace
Equation (\ref{Bohr-Sommerfeld}) gives the Bohr-Sommerfeld conditions
discussed in the \linebreak 
introduction. The action integral is independent of the
parametrization of $\gamma $ within its orientation class. However, the
change of orientation of $\gamma $ would lead to the change from $n$ to $-n$. Therefore, the Bohr-Sommerfeld condition (\ref{Bohr-Sommerfeld}) depends only on the image of $\gamma $. In the following, we shall refer to the image an integral curve $\gamma $ of $X_{H}$ that satisfies equation (\ref{Bohr-Sommerfeld}) as a \emph{Bohr-Sommerfeld torus}. The integer $n$ on the right hand side of equation (\ref{Bohr-Sommerfeld}) is called the \emph{quantum number}
of the corresponding Bohr-Sommerfeld torus. Since integral curves of $X_{H}$
preserve the Hamiltonian $H$, we may rewrite equation (\ref{Bohr-Sommerfeld}) in the form 
\begin{equation}
I|_{C(e)}=\ttfrac{1}{2\pi} \int_{C(e)}p~\mathrm{d}\alpha =n\hbar ,  
\label{Bohr-Sommerfeld1}
\end{equation}
where $C(e)$ is a connected component of the energy level $H^{-1}(e)$. Thus,
Bohr-Sommerfeld conditions (\ref{Bohr-Sommerfeld}) impose conditions on the
energy. The set of values of the energy allowed by Bohr-Sommerfeld
conditions is interpreted as the \emph{quantum energy spectrum} of the system. \medskip 

\noindent From the discussion preceding theorem \ref{BS} it follows that a 
section $\sigma $ of the prequantum line bundle, which is covariantly
constant along $D$, has support contained in the union of Bohr-Sommerfeld tori and the energy
levels $H^{-1}(0)$ and $H^{-1}(2)$. Since $H^{-1}(0)$ is a critical point, the restriction of $\sigma $ to 
$H^{-1}(0)$ is the value of $\sigma $ at $H^{-1}(0)$ which is not restricted
by the condition that $\sigma $ is covariantly constant along $D$. So we
may allow the value $e=0$ in equation (\ref{Bohr-Sommerfeld1}). On the one hand, we consider 
$H^{-1}(0)$ as a (singular) Bohr-Sommerfeld torus corresponding to the quantum
number $n=0$. On the other hand, we assume that the singular level set $H^{-1}(2)$ is \emph{not} a (singular) Bohr-Sommerfeld torus. \medskip 

Since a section of theprequantization line bundle that is covariantly constant along $D$ has its support in the union of Bohr-Sommerfeld tori, which has empty interior, such sections can be smooth only in the sense of distributions. Therefore, we adopt the following definition. \medskip

\noindent \textbf{Definition 4.2} \textit{A quantum state of the mathematical pendulum 
is a section $\sigma $ of the
prequantization line bundle $\rho $, whose support lies in the union of
Bohr-Sommerfeld tori such that for each Bohr-Sommerfeld torus $C$ the
restriction $\sigma |_C$ of $\sigma $ to $C$ is a smooth covariantly
constant section of $\rho |_C$.} \medskip 

\noindent Let $\mathfrak{H}$ be the space of quantum states of the mathematical
pendulum. For each Bohr-Sommerfeld torus $C$, we choose a non-vanishing
smooth covariantly constant section $\sigma $ of $L_{\mid C}$. The family $\{\sigma _{\mid C}\}$ is 
a basis of $\mathfrak{H}$, which we shall refer to as a
Bohr-Sommerfeld basis. Give $\mathfrak{H}$ a hermitian scalar product $(\cdot \mid
\cdot )$ so that the Bohr-Sommerfeld basis $\{ {\sigma}_{\mid C}\}$ is orthonormal. 
Thus, we have obtained a vector space structure on the space of states of
the mathematical pendulum. Note that this structure is not uniquely
determined by the geometry of the classical phase space. We have the freedom
of multiplying each basis vector $\sigma _{\mid C}$ by a nonzero complex number. \medskip 

\noindent \textbf{Definition 4.3} \textit{A function $f\in C^{\infty }(T^{\ast }S^{1})$ is Bohr-Sommerfeld quantizable if it is constant on Bohr-Sommerfeld tori. Bohr-Sommerfeld
quantization assigns to a quantizable function $f$ a linear operator 
$\boldsymbol{Q}_{f}$ on $\mathfrak{H}$ such that, for each Bohr-Sommerfeld
torus $C$} 
\begin{equation}
\boldsymbol{Q}_{f}\sigma _{\mid C}=f_{\mid C}\, \, \sigma _{\mid C}.  \label{BSquantization}
\end{equation}

Observe that the operators $\boldsymbol{Q}_{f}$ corresponding to
Bohr-Sommerfeld quantizable functions $f$ are diagonal in the
Bohr-Sommerfeld basis. Since the Bohr-Sommerfeld tori are closed and
mutually disjoint, for any function $C\mapsto \lambda _{C}$ on the
collection of Bohr-Sommerfeld tori, there exists a function $f\in C^{\infty}(T^{\ast }S^{1})$ such that $f_{\mid C}=\lambda _{C}$. Thus, each basis vector $\sigma _{\mid C}$ is an eigenvector of the operator $\boldsymbol{Q}_{f}$ corresponding to an eigenvalue $\lambda _{C}$.

\section{Quantization away from the singularity}\label{section-5}
%%%%%%%%%%%%%%%%%%%%%

\subsection{Structure of the Bohr-Sommerfeld basis}\label{section-5.1}
%%%%%%%%%%%%%%%%%%%

We now study of the structure of the Bohr-Sommerfeld basis $\{\sigma |_C \}.$ \medskip 

The energy level $H^{-1}(2)$ divides $T^{\ast }S^{1}$ into three open
subsets: $P_{0} = \{(p,\alpha )\mid H(p,\alpha )<2\} $ and $P_{\mp} = \{(p,\alpha )\mid H(p,\alpha )>2\text{ and }\mp p>0\} $. Let $\mathfrak{H}_{0}$ be the subspace of $\mathfrak{H}$ consisting of
sections with support in $P_{0}$. Similarly, let $\mathfrak{H}_{\mp}$ be the subspaces of 
$\mathfrak{H}$ consisting of sections with support in $P_{\mp }$. Then 
$\mathfrak{H}=\mathfrak{H}_{0}\oplus \mathfrak{H}_{-}\oplus \mathfrak{H}_{+}$ and 
${\{ {\sigma }|_C  \}}_{0, \pm } = \{ \sigma |_C \, \, \setrule \, \, C\subseteq P_{0,\pm} \} $ are bases of $\mathfrak{H}_{0}$, ${\mathfrak{H}}_{+}$, and ${\mathfrak{H}}_{-}$, respectively. \medskip

A Bohr-Sommerfeld torus in $P_0$ can be labelled by its quantum number $n = 0, \ldots , N$, where $N$ is the largest nonnegative integer such that $N\hbar < I(2)$. Thus 
${\{ \sigma |_C \} }_0 = \{ {\sigma }^0_0, \ldots , {\sigma }^0_N \}$, where the subscript 
$n =0, \ldots , N$ is the quantum number of the state ${\sigma }^0_n$ and the superscript $0$ reminds us that ${\sigma }^0_n$ lies in ${\mathfrak{H}}_0$. 
Similarly a Bohr-Sommerfeld torus in $P_{\pm}$ can be 
labeled by it quantum number $m \ge M$, where $2M$ is the smallest even nonnegative integer 
greater than or equal to $N+1$. In other words, 
\begin{equation}
M = \min \{ m \in {\Z }_{>0} \setrule 2m \ge N+1 \} = 
\left\{ \begin{array}{rl} 
\onehalf (N+2), & \mbox{if $N$ is even} \\
\rowspace \onehalf (N+1), & \mbox{if $N$ is odd.} \end{array} \right. 
\label{eq-news5one}
\end{equation}
Hence the Bohr-Sommerfeld basis of ${\mathfrak{H}}_{\pm }$ is ${\{ \sigma |_C \} }_{\pm } = 
\{ {\sigma }^{\pm }_m, \, \, m \ge M \} $, where the subscript $m$ is the quantum number of 
the state ${\sigma }^{\pm }_m$ in ${\mathfrak{H}}_{\pm }$. Since $N < 2M$ the basis 
$\{ {\sigma }_C \} $ of $\mathfrak{H}$ has the lattice structure 
\begin{equation}
{\sigma }^0_0 \leftrightarrows {\sigma }^0_1 \leftrightarrows \cdots 
\leftrightarrows {\sigma }^0_N 
\begin{array}{l}
\nearrow \! \! \! \swarrow \rule{0pt}{10pt}\raisebox{5pt}{${\sigma }^{+}_M \leftrightarrows 
{\sigma }^{+}_{M+1} \leftrightarrows \cdots $} \\
\searrow \! \! \! \nwarrow \raisebox{-4pt}{${\sigma }^{-}_M \leftrightarrows
{\sigma }^{-}_{M+1} \leftrightarrows \cdots $}
\end{array}
\label{eq-news5two}
\end{equation}

The structure of the Bohr-Sommerfeld set can be used to
study the energy spectrum of the mathematical pendulum, which  consists of
values of $e_{n}$ such that a connected component $C(e_{n})$ of the energy
level $H^{-1}(e_{n})$ satisfies Bohr-Sommerfeld conditions%
\begin{displaymath}
\frac{1}{2\pi }\int_{C(e_{n})}pd\alpha =n\hbar 
\end{displaymath}
for some integer $n$. See the discussion following equation (\ref{nabla}). The part
of energy spectrum contained in the interval interval $[0,2]$ is simple, and
it can be obtained by solving for $e_{n}$ equation 
\begin{displaymath}
n\hbar =\frac{4}{\pi }e_{n}\int_{0}^{\pi /2}\frac{\cos ^{2}\varphi }{\sqrt{1-%
\frac{e_{n}}{2}\sin ^{2}\varphi }}d\varphi ,
\end{displaymath}
where $0\leq n\leq N$, and $N$ is the largest positive integer such that 
$e_{N}<2$. We have assumed that $e=2$ is not in the energy spectrum of the
mathematical pendulum. For $e>2$, The part of the energy spectrum contained
in the half line $[2,\infty )$ can be obtained by solving for $e_{m}$ the equation 
\begin{displaymath}
m\hbar =\frac{1}{\pi }\int_{-\pi }^{\pi }\sqrt{2\left( e_{m}-(1-\cos \alpha
\right) }d\alpha ,
\end{displaymath}
where $2m\geq N+1\ $ensures that $e_{m}>2$. In this range, each eigenspace
is 2-dimensional. 

\subsection{Transitions between quantum states}\label{section-5.2}
%%%%%%%%%%%%%%%%%%%

In this subsection we discuss the transitions between quantum states given by the horizontal 
arrows in diagram (\ref{eq-news5two}). The transitions from ${\sigma }^0_N$ to 
${\sigma }^{\pm }_M$ given by slanted arrows in diagram (\ref{eq-news5two}) 
involve crossing the energy level $2$, where the action $I$ is continuous but not 
differentiable. This requires understanding the ${\Z }_2$ symmetry of the mathematical 
pendulum, which will be treated in the next section. \medskip  

In diagram (\ref{eq-news5two}) transitions involving the right pointing horizontal arrows correspond to the action of an operator $\mathbf{b}$ on $\mathfrak{H}$ such that 
\begin{equation}
\mathbf{b}\, {\sigma }^0_n = {\sigma }^0_{n+1}, \, \, \mbox{for $n =0, 1, \ldots , N-1$}
\label{eq-news5three}
\end{equation}
and 
\begin{equation}
\mathbf{b}\, {\sigma }^{\pm }_m = {\sigma }^{\pm }_{m+1}, \, \, \mbox{for $m = M, M+1, \ldots $.}
\label{eq-news5four}
\end{equation}
We refer to $\mathbf{b}$ as the \emph{raising} operator on $\mathfrak{H}$. 
Transitions involving the left pointing horizontal arrows give rise to the 
\emph{lowering} operator $\mathbf{a}$ such that 
\begin{equation}
\mathbf{a}\, {\sigma }^0_n  = {\sigma }^0_{n-1}, \, \, \mbox{for $n=1, \ldots , N$} 
\label{eq-news5twelve} 
\end{equation}
and
\begin{equation}
\mathbf{a}\, {\sigma }^{\pm }_m  = {\sigma }^{\pm }_{m-1}, \, \, \mbox{for $m = M+1, \ldots $.}
\label{eq-news5thirteen}
\end{equation}
Since ${\sigma }^0_0$ is the lowest point in the lattice, we require that 
\begin{equation}
\mathbf{a}\, {\sigma }^0_0 = 0.
\label{eq-news5fourteen}
\end{equation}
Thus in diagram (\ref{eq-news5two}) the lowering operator $\mathbf{a}$ corresponds to left pointing horizontal arrows, while the raising operator $\mathbf{b}$ corresponds to right pointing horizontal arrows.\footnote{For small positive quantum numbers, when the value of $e$ is slightly above zero, the mathematical pendulum is well approximated by the 
$1$-dimensional harmonic oscillator. In this situation the shifting operators 
$\mathbf{a}$ and $\mathbf{b}$ are well appoximated by the classical raising and 
lowering operators of the harmonic oscillator. We do not discuss the complex analytic nature of this approximation.} In order to interpret the slanted arrows in the diagram we need to discuss the ${\Z }_2$ symmetry of the mathematical pendulum. The shifting operators $\mathbf{a}$ and 
$\mathbf{b}$ will be constructed by lifting the shifting operator in the ${\Z }_2$-reduced 
quantum system. In the appendix we construct this shifting operator using 
geometric quantization.   \medskip 

In order to identify the function whose quantization might lead to the operator $\mathbf{b}$ we extend Dirac's quantization rule 
\begin{equation}
[{\mathbf{Q}}_{f_1}, {\mathbf{Q}}_{f_2}] = i\hbar \, {\mathbf{Q}}_{\{ f_1, f_2 \} } 
\label{eq-news5five}
\end{equation}
to complex valued functions. Action angle coordinates $(I, \vartheta )$ on $P_0 \cup P_{+} \cup P_{-}$ restrict to action angle coordinates $(I_0, {\vartheta }_0)$ on $(P_0, {\omega }|_{P_0})$  and $(I_{\pm}, {\vartheta }_{\pm})$ on $(P_{\pm}, {\omega }|_{P_{\pm}})$, respectively. The latter action angle coordinates have been computed in section 2.2. They satisfy the Poisson bracket relations 
$\{ I_0, {\vartheta }_0 \} = -1$ and $\{ I_{\pm }, {\vartheta }_{\pm } \} = -1$ on $P_0$ and 
$P_{\pm }$, respectively. Therefore
\begin{equation}
\{ I_0, {\mathrm{e}}^{i{\vartheta }_0 } \} = -i{\mathrm{e}}^{i{\vartheta }_0} \, \, \, 
\mathrm{and} \, \, \, 
\{ I_{\pm}, {\mathrm{e}}^{i{\vartheta }_{\pm} } \} = -i{\mathrm{e}}^{i{\vartheta }_{\pm}} . 
\label{eq-news5six}
\end{equation}
If we introduce the quantum operator ${\mathbf{Q}}_{{\mathrm{e}}^{i \vartheta }}$ on 
$\mathfrak{H}$, equations (\ref{eq-news5five}) and (\ref{eq-news5six}) imply 
\begin{equation}
[{\mathbf{Q}}_{I_0}, {\mathbf{Q}}_{{\mathrm{e}}^{i{\vartheta }_0}}] = 
-i\hbar \, {\mathbf{Q}}_{{\mathrm{e}}^{i{\vartheta }_0}} \, \, \, 
\mathrm{and} \, \, \, 
[{\mathbf{Q}}_{I_{\pm }}, {\mathbf{Q}}_{{\mathrm{e}}^{i{\vartheta }_{\pm}}}] = 
-i\hbar \, {\mathbf{Q}}_{{\mathrm{e}}^{i{\vartheta }_{\pm }}} . 
\label{eq-news5seven}
\end{equation}
On the other hand, equation (\ref{BSquantization}) yields 
\begin{equation}
{\mathbf{Q}}_{I_0}{\sigma }^0_n = n\hbar \, {\sigma }^{0}_n, \, \, \mbox{for $n =0, 1, \ldots , N$} 
\label{eq-news5eight}
\end{equation}
and 
\begin{equation}
{\mathbf{Q}}_{I_{\pm}}{\sigma }^{\pm}_m = m\hbar \, {\sigma }^{\pm}_m, \, \, 
\mbox{for $m =M, M+1, \ldots $} 
\label{eq-news5nine}
\end{equation}
Equations (\ref{eq-news5three}) and (\ref{eq-news5four}) imply that the operator $\mathbf{b}$ and the quantized actions satisfy the commutation relations 
\begin{equation}
[{\mathbf{Q}}_{I_0}, \mathbf{b} ] {\sigma }^0_n  = -i\hbar \, \mathbf{b}\, {\sigma }^0_n \, \, \, 
\mathrm{and} \, \, \, 
[{\mathbf{Q}}_{I_{\pm }}, \mathbf{b} ] {\sigma }^{\pm}_m  = -i\hbar \, \mathbf{b}\, 
{\sigma }^{\pm }_m 
\label{eq-news5ten}
\end{equation}
for $n =0 ,1, \ldots , N-1$ and $m = M, M+1, \ldots $, respectively. Comparing equations 
(\ref{eq-news5seven}) and (\ref{eq-news5ten}) shows that the raising operator $\mathbf{b}$ defined in equations (\ref{eq-news5three}) and (\ref{eq-news5four}) satisfy the same commutation relations as the operator ${\mathbf{Q}}_{{\mathrm{e}}^{i\vartheta }}$. \medskip  

It is of interest to see to what degree the raising and the lowering
operators on $\mathfrak{H}$ correspond to quantization of classical
functions on $P=P_{0}\cup P_{-}\cup P_{+}$. Their restrictions to the
subspaces $\mathfrak{H}_{0},$ ${\mathfrak{H}}_{-}$ and ${\mathfrak{H}}_{+}$ of 
$\mathfrak{H}$ with supports in ${\mathfrak{H}}_{0}$, ${\mathfrak{H}}_{-}$ and 
${\mathfrak{H}}_{+}$, respectively, can be related to quantization of the functions 
${\mathrm{e}}^{\pm i\vartheta _{0}}$, ${\mathrm{e}}^{\pm i{\vartheta }_{-}},$ and 
${\mathrm{e}}^{\pm i{\vartheta }_{-}}$ on $P_{0},$ $P_{+}$ and $P_{-}$, respectively. \medskip 

Choose the basic sections of basic sections ${\{ {\sigma }^0_n \} }^N_{n=1}$ in 
${\mathfrak{H}}_{0}$ so that  
${\mathrm{e}}^{i{\vartheta }_{0}}{\sigma}^{0}_{n} = 
\mathbf{a}\, {\sigma }^{0}_{n+1}$ for all $1 \le n \le N$. Because 
${\mathbf{b}}_{0} \comp {\mathbf{a}}_{0} \, {\sigma }^{0}_n = 
{\sigma }^{0}_n$ for all $1 \le n \le N$, it follows that the restriction ${\mathbf{b}}_0$ of the raising operator $\mathbf{b}$ to ${\mathfrak{H}}_{0}$ may be interpreted as the operator 
${\mathbf{Q}}_{{\mathrm{e}}^{i {\vartheta }_{0}}}$, which sends ${\sigma }^0_n$ into 
${\mathrm{e}}^{i {\vartheta }_{0}}{\sigma }^0_n$ for $1 \le n \le N$. Similarly, the restriction 
${\mathbf{b}}_{\mp}$ of the raising operator $\mathbf{b}$ to ${\mathfrak{H}}_{-}$ may be
interpreted as the operator ${\mathbf{Q}}_{{\mathrm{e}}^{i {\vartheta }_{-}}}$
of multiplication by ${\mathrm{e}}^{i{\vartheta }_{0, \mp }}$. In other words, 
${\mathbf{b}}_{-}= {\mathbf{Q}}_{{\mathrm{e}}^{i {\vartheta }_{-}}}$.
Similarly, ${\mathbf{b}}_{+}= {\mathbf{Q}}_{{\mathrm{e}}^{i {\vartheta }_{+}}}$, 
is the restriction of the raising operator $\mathbf{b}$ to 
${\mathfrak{H}}_{+.}$ On the other hand, the restriction ${\mathbf{b}}_{0}$ of the
raising operator $\mathbf{b}$ to ${\mathfrak{H}}_{0}$ agrees with 
${\mathrm{e}}^{i {\vartheta }_{0}}$ except on ${\sigma }^{0}_{n}$ because 
$\mathbf{b}{\sigma }^0_{n}$ is not in ${\mathfrak{H}}_{0}$. In the same way, the lowering
operators ${\mathbf{a}}_{0},$ ${\mathbf{a}}_{-}$, and ${\mathbf{a}}_{+}$ 
may be interpreted in terms of the multiplication operators 
${\mathbf{Q}}_{{\mathrm{e}}^{i{\vartheta }_{0}}}$, 
${\mathbf{Q}}_{{\mathrm{e}}^{i {\vartheta }_{-}}}$, and 
${\mathbf{Q}}_{{\mathrm{e}}^{i {\vartheta } _{+}}}$, respectively.

\section{The ${\Z }_2$-symmetry} \label{section-6}
%%%%%%%%%%%%%%%%%%%

The mathematical pendulum $(H, T^{\ast }S^1, \omega )$ has a ${\Z }_2$-symmetry generated by 
\begin{equation}
\zeta : T^{\ast }S^1 \rightarrow T^{\ast }S^1: (p, \alpha ) \mapsto (-p, -\alpha ), 
\label{eq-nwzero}
\end{equation}
because the Hamiltonian $H$ (\ref{eq-s1one}) and the symplectic form $\omega = 
\dee p \wedge \dee \alpha $ are invariant. In more detail, 
for every $(p, \alpha ) \in T^{\ast }S^1$ we have 
$({\zeta }^{\ast }H)(p, \alpha ) = H(-p,-\alpha ) = H(p, \alpha )$. So $H$ is invariant. Also 
the $1$-form $\theta = p \dee \alpha $ is invariant, because $({\zeta }^{\ast }\theta )(p,\alpha ) 
= (-p) \dee (-\alpha ) = \theta (p, \alpha ) $. This implies that the $2$-form $\omega $ is invariant, 
because ${\zeta }^{\ast }\omega = {\zeta }^{\ast }(\dee \theta ) =  
\dee \, ({\zeta }^{\ast }\theta )  = \omega $. \medskip 

The quantized mathematical pendulum $(H, T^{\ast }S^1, \omega )$ has quantum line bundle 
\begin{equation}
\rho : L = \C \times T^{\ast }S^1 \rightarrow T^{\ast }S^1: \big( z, (p, \alpha ) \big) \mapsto 
(p, \alpha )
\label{eq-nws6zero}
\end{equation}
with covariant derivative $\nabla = -i\hbar \theta \otimes {\lambda }_0$, where 
${\lambda }_0: T^{\ast }S^1 \rightarrow L: (p, \alpha ) \mapsto \big( 1, (p, \alpha ) \big) $ is 
a section, which trivializes the bundle $\rho $. The bundle space $L$ has a ${\Z }_2$-symmetry generated by the mapping 
\begin{equation}
\mu : L \rightarrow L : \big( z, (p, \alpha ) \big) \mapsto \big( z, (-p, -\alpha ) \big) . 
\label{eq-nws6one}
\end{equation}
The mapping $\mu $ covers the ${\Z }_2$-symmetry of the mathematical pendulum generated 
by $\zeta $ (\ref{eq-nwzero}), since $\rho \comp \mu = \zeta \comp \rho$.

\subsection{The ${\Z }_2$-symmetric quantum system}\label{section-6.1}
%%%%%%%%%%%%%%%%%%%%%

Consider the ${\Z }_2$-symmetry on $L$ generated by the mapping $\mu $ 
(\ref{eq-nws6one}). The ${\Z }_2$-symmetric quantized system is the mathematical 
pendulum with ${\Z}_2$-symmetry generated by $\zeta $ (\ref{eq-nwzero}) and 
quantum line bundle $\rho $ (\ref{eq-nws6zero}) having the ${\Z}_2$-symmetry 
generated by $\mu $ (\ref{eq-nws6one}). \medskip %

Let $\Gamma (\rho )$ be the vector space of smooth sections of the bundle $\rho $. \medskip   

\noindent \textbf{Lemma 6.1.1} \textit{The ${\Z }_2$ action on $L$ generated by the 
mapping $\mu $ $\mathrm{(\ref{eq-nws6one})}$ induces a ${\Z }_2$-action on 
$\Gamma (\rho )$ generated by the linear map}
\begin{equation}
{\mu }^{\ast }: \Gamma (\rho ) \rightarrow \Gamma (\rho ): \sigma \mapsto {\mu }^{\ast }(\sigma ) 
= {\mu }^{-1} \comp \sigma \comp \zeta . 
\label{eq-nws6onestar}
\end{equation}

\begin{proof} Let $\sigma : T^{\ast }S^1 \rightarrow L: (p, \alpha ) \mapsto \big( z(p,\alpha ), (p, \alpha ) \big) $ be a smooth section of the bundle $\rho $. Then 
\begin{equation}
({\mu }^{\ast }\sigma )(p, \alpha ) = 
{\mu }^{-1}\big( z\big( \zeta (p, \alpha ) \big), \zeta (p, \alpha ) \big) = 
\big( z(\zeta (p, \alpha )), (p, \alpha ) \big) 
\label{eq-nws6onedot}
\end{equation} 
is a section of the bundle $\rho $. So 
\begin{displaymath}
{\mu }^{\ast }\big( {\mu }^{\ast }(\sigma )\big) (p, \alpha ) = 
\big( z({\zeta }^2(p,\alpha )) , (p,\alpha ) \big) = \sigma (p, \alpha ), 
\end{displaymath}
that is, $({\mu }^{\ast})^2 \sigma = \sigma $. \end{proof} \medskip 

The mapping $\zeta :T^{\ast } S^1 \rightarrow T^{\ast }S^1$ (\ref{eq-nwzero}) acts on the set 
of all Bohr-Sommerfeld tori. Hence it induces an operator 
\begin{equation}
\mathbf{P}: \mathfrak{H} \rightarrow \mathfrak{H}:\sigma |_C \mapsto {\mu }^{\ast }(\sigma |_C) ,
\label{eq-s6nthree}
\end{equation}
which we call the \emph{parity} operator. Since the mapping $\mu $ (\ref{eq-nws6one}) 
generates a ${\Z }_2$-action on L, it follows that ${\mu }^{\ast }$ (\ref{eq-nws6onestar}) generates a representation of ${\Z }_2$ on $\mathfrak{H}$. From the fact that the Hamiltonian $H$ (\ref{eq-s1one}) is invariant under the ${\Z }_2$-action on $T^{\ast }S^1$, we get  $[{\mathbf{Q}}_H, \mathbf{P}] = 0$. 
Moreover, \medskip 

\noindent \textbf{Lemma 6.1.2} \label{restrictedparityops}The maps  
\begin{displaymath}
{\mathbf{P}}_{0, \pm} = \mathbf{P}|_{{\mathfrak{H}}_{0,\pm }}: {\mathfrak{H}}_{0,\pm } 
\rightarrow {\mathfrak{H}}_{0, \mp } 
\end{displaymath}
are bijective involutions, which implies ${\mathbf{P}}^{-1}_0 = {\mathbf{P}}_0$ 
and ${\mathbf{P}}^{-1}_{\pm } = 
{\mathbf{P}}_{\mp }$.  

\begin{proof} Recall that if $\sigma \in \Gamma (\rho )$, then $\sigma = f {\lambda }_0$ for 
some smooth function $f$ on $T^{\ast }S^1$. The support $\mathrm{supp}\, \sigma $ of $\sigma $ is the support $\mathrm{supp}\, f$ of the function $f$, which is $\{ (p, \alpha ) \in T^{\ast }S^1 \setrule 
f(p,\alpha ) \ne 0 \} $. Suppose that ${\sigma }|_{C_{0,\pm}(e)} = 
z_{0,\pm}(p,\alpha ) {\lambda }_0(p, \alpha )
\in {\mathfrak{H}}_{0,\pm }$. Then $\mathrm{supp}\, {\sigma }|_{C_{0,\pm}(e)} = C_{0,\pm}(e)$. 
From (\ref{eq-nws6onedot}) we get $\big( {\mathbf{P}}_{0, \pm }(\sigma |_{C_{0,\pm}(e)}) \big) (p, \alpha ) = z_{0,\pm}(\zeta (p, \alpha )) {\lambda }_0(p,\alpha )$. So we obtain 
\begin{displaymath}
(p,\alpha ) \in 
\mathrm{supp}\, {\mathbf{P}}_{0,\pm }(\sigma |_{C_{0,\pm }(e)}) \, \, \, 
\mathrm{if \, and\, only\, if} \, \, \,  
\zeta (p, \alpha ) \in \mathrm{supp}\, z_{0,\pm } = C_{0,\pm}(e)
\end{displaymath}
if and only if 
$(p, \alpha ) \in \zeta \big( C_{0,\pm}(e) \big) = C_{0, \mp}(e)$. So 
${\mathbf{P}}_{0, \pm }(\sigma |_{C_{0,\pm}(e)}) \in {\mathfrak{H}}_{0,\mp}$. Thus 
${\mathbf{P}}_{0,\pm }$ maps ${\mathfrak{H}}_{0,\pm }$ onto ${\mathfrak{H}}_{0,\mp}$. 
From 
\begin{align}
\big( {\mathbf{P}}^2_{0, \pm}\sigma |_{C_{0,\pm}(e)} \big) (p, \alpha ) & = 
z_{0,\pm }\big( {\zeta }^2(p, \alpha ) \big) {\lambda }_0(p,\alpha ) \notag \\
& = z_{0,\pm }(p, \alpha ) {\lambda }_0(p,\alpha ) = \sigma |_{C_{0,\pm}(e)}(p, \alpha ) \notag 
\end{align}
we get ${\mathbf{P}}^2_{0,\pm} = {\mathrm{id}}|_{{\mathfrak{H}}_{0,\pm}}$. The operator 
${\mathbf{P}}_{0,\pm }$ is injective, for if ${\mathbf{P}}_{0,\pm }( \sigma |_{C_{0,\pm}(e)} )$  
$= {\mathbf{P}}_{0,\pm }( {\sigma }'|_{C_{0,\pm}(e)} )$, then 
\begin{displaymath}
\sigma |_{C_{0,\pm}(e)} = {\mathbf{P}}^2_{0,\pm }( \sigma |_{C_{0,\pm}(e)} )  = 
{\mathbf{P}}^2_{0,\pm }( {\sigma }' |_{C_{0,\pm}(e)} ) = {\sigma }' |_{C_{0,\pm}(e)}. 
\end{displaymath}
Thus the operator ${\mathbf{P}}_{0,\pm }$ is bijective. Since ${\mathbf{P}}^2_{0,\pm } = 
{\mathrm{id}}_{{\mathfrak{H}}_{\pm}}$, it follows that ${\mathbf{P}}^{-1}_0 = {\mathbf{P}}_0$ and 
${\mathbf{P}}^{-1}_{\pm } = {\mathbf{P}}_{\mp}$. \end{proof} \medskip 

We say that a quantum state ${\sigma }|_C$ of the ${\Z}_2$-quantized mathematical pendulum 
$(H, T^{\ast }S^1, \omega )$ with quantum line bundle $\rho $ is \emph{even} if it is covariantly constant even section $\sigma |_C$ of $\rho $, that is, $\mathbf{P}(\sigma |_C) = \sigma |_C$. Let
${\mathfrak{H}}^{\mathrm{even}}$ be the vector space spanned by the even quantum states. 
${\sigma }|_C$ is an \emph{odd} quantum state if it is covariantly constant odd section 
$\sigma |_C$ of $\rho $, that is, $\mathbf{P}(\sigma |_C) = - \sigma |_C$. Let
${\mathfrak{H}}^{\mathrm{odd}}$ be the vector space spanned by the odd quantum states. 
\medskip  

\noindent \textbf{Claim 6.1.3} \textit{We have}  
\begin{equation}
\mathfrak{H} = {\mathfrak{H}}^{\mathrm{even}} \oplus {\mathfrak{H}}^{\mathrm{odd}}. 
\label{eq-sn6four}
\end{equation}

\begin{proof} The proof of this claim is standard, but we include it for completeness. 
Suppose that the section $\sigma |_C$ is in $\mathfrak{H}$. Then 
$\mathbf{P}(\sigma |_C) \in \mathfrak{H}$. So the sections 
${\sigma}_{\mathrm{even}}|_C = \onehalf \big( \sigma |_C + \mathbf{P}(\sigma |_C) \big) $ and 
${\sigma}_{\mathrm{odd}}|_C = \onehalf \big( \sigma |_C - \mathbf{P}(\sigma |_C) \big) $ lie 
in $\mathfrak{H}$ and $\sigma |_C = {\sigma }_{\mathrm{even}}|C + 
{\sigma }_{\mathrm{odd}}|_C$. Now ${\sigma }_{\mathrm{even}}|_C \in {\mathfrak{H}}^{\mathrm{even}}$,
because it is covariantly constant, has support $C$, and is an even section, since 
\begin{displaymath}
\mathbf{P}({\sigma }_{\mathrm{even}}|_C) = \onehalf \big( \mathbf{P}(\sigma |_C) + 
{\mathbf{P}}^2(\sigma |_C) \big) = \onehalf \big(\sigma |_C + \mathbf{P}(\sigma |_C) \big) 
= {\sigma }_{\mathrm{even}}|_C. 
\end{displaymath}
Similarly, ${\sigma }_{\mathrm{odd}}|_C \in {\mathfrak{H}}^{\mathrm{odd}}$. Thus 
$\mathfrak{H} = {\mathfrak{H}}^{\mathrm{even}} + {\mathfrak{H}}^{\mathrm{odd}}$. The preceding 
sum is direct since ${\mathfrak{H}}^{\mathrm{even}} \cap {\mathfrak{H}}^{\mathrm{odd}} = 
\{ 0 \} $. For if $\sigma |_C \in {\mathfrak{H}}^{\mathrm{even}} \cap 
{\mathfrak{H}}^{\mathrm{odd}}$, 
then $\sigma |_C = \mathbf{P}(\sigma |_C) = - \sigma |_C$, which implies $\sigma |_C =0$. 
\end{proof} \medskip 

\noindent \textbf{Theorem 6.1.4} \textit{Let ${\sigma }|_C$ be a nonzero even or odd quantum state in ${\mathfrak{H}}^{\mathrm{even}}\cap {\mathfrak{H}}_0$ or 
${\mathfrak{H}}^{\mathrm{odd}} \cap {\mathfrak{H}}_0$, respectively. 
Then its quantum number is even or odd, respectively.}  \medskip  

\begin{proof} Write $\sigma (p, \alpha ) = z(p, \alpha ){\lambda }_0$ for $(p, \alpha )  
\in T^{\times }S^1= T^{\ast }S^1\setminus \{ (0,0), (0, \pi ) \}$, 
where $z(p, \alpha ) = \varepsilon z(-p,-\alpha )$ with $\varepsilon =1$ if 
$\sigma $ is even and $\varepsilon = -1$ if $\sigma $ is odd. Note that $z|_C$ is nowhere 
vanishing, for if it vanished at some point in $C$ then it would be identically zero on $C$, 
since $\sigma |_C$ is covariantly constant. But this contradicts our hypothesis. 
Because $\sigma |_C \in {\mathfrak{H}}_0$ by hypothesis, it follows that 
$C = C(e) = (H^{\times })^{-1}(e)$, where $H^{\times } = H|T^{\times }S^1$ and $0 < e < 2$. \medskip 

Let $\gamma $ be a closed curve in $T^{\times }S^1$, which parametrizes $C$. The 
image of $\gamma $ is $\{ (p, \alpha ) \in T^{\times }S^1 \setrule \onehalf p^2 - \cos \alpha +1 =e \} $, 
which is clearly invariant under the ${\Z }_2$-action generated by $\zeta $ (\ref{eq-nwzero}). 
Setting $p=0$ we get $1 -\cos \alpha = e$, which has a solution ${\alpha }^{+} 
\in [0, \pi ]$ and another solution ${\alpha }^{-} = - {\alpha }^{+}$. Thus we can 
parametrize $\gamma $ by $\alpha $, say $\gamma (\alpha ) = 
\big( p(\alpha ) , \alpha \big) $ for $\alpha \in [-{\alpha }^{+}, {\alpha }^{+}]$. Let 
\begin{displaymath}
{\rho }^{\times } : L^{\times } = \C \times T^{\times }S^1 \rightarrow T^{\times }S^1: 
\big( z, (p, \alpha ) \big) \mapsto (p, \alpha )
\end{displaymath}
be the quantum line bundle with covariant derivative ${\nabla }^{\times }= 
\nabla |\Gamma ({\rho }^{\times })$. Let 
$\widetilde{\gamma }$ be the horizontal lift of $\gamma $. Parametrize the image of 
$\widetilde{\gamma }$ by $\alpha $, namely, $\widetilde{\gamma }(\alpha ) = 
\big( z(p(\alpha ), \alpha ), (p(\alpha ), \alpha ) \big) $ for $\alpha \in [-{\alpha }^{+}, {\alpha }^{+}]$. So $\widetilde{\gamma }$ is a section of $L^{\times }$. \medskip 

Suppose that for some ${\alpha }_0 \in [-{\alpha }^{+}, {\alpha }^{+}] $ we have 
$z\big( \gamma (-{\alpha }_0) \big) = \varepsilon z \big( \gamma ({\alpha }_0) \big)$, where 
$\varepsilon = \pm 1$. Then 
$z\big( \gamma (-{\alpha }_0) \big) \ne 0$, since $\sigma |_C$ is nowhere zero. Because  
$\widetilde{\gamma }$ is the horizontal lift of $\gamma $, using the covariant derivative 
${\nabla }^{\times }$, integrating equation (\ref{paralleltransport}) we get 
\begin{align}
-i\frac{\hbar }{2\pi } \ln \frac{z(\gamma ({\alpha }_0))}{z(\gamma (-{\alpha }_0))} & = 
\frac{2}{2\pi } \int^{{\alpha }_0}_{-{\alpha }_0} \sqrt{2(e-(1-\cos \alpha ))} \, \dee \alpha \notag \\
& = \left\{ \begin{array}{rl} 
2m\hbar , & \mbox{if $\varepsilon =1$} \\
(2m-1)\hbar , & \mbox{if $\varepsilon = -1$} 
\end{array} \right. 
\label{eq-news6five}%
\end{align}
for some positive integer $m$. The last equality in (\ref{eq-news6five}) follows because 
the hypothesis $z\big( \gamma (-{\alpha }_0) \big) = 
\varepsilon z\big( \gamma ({\alpha }_0) \big) $ 
implies that 
\begin{displaymath}
\ln \frac{z(\gamma ({\alpha }_0))}{z(\gamma (-{\alpha }_0))} = 
\ln \varepsilon = \left\{ \begin{array}{rl} 
2m\pi i , & \mbox{if $\varepsilon =1$} \\
(2m-1)\pi i , & \mbox{if $\varepsilon = -1$}
\end{array} \right. 
\end{displaymath}
for some integer $m$. Since $\sqrt{2(e-(1-\cos \alpha ))} >0$ on $(-{\alpha }_0, {\alpha }_0)$, 
we see that $m >0$. From the fact that the Bohr-Sommerfeld torus $C = C(e) = 
(H^{\times })^{-1}(e)$ has quantum number $n$, that is, 
\begin{displaymath}
\frac{2}{2\pi } \int^{{\alpha }^{+}}_{-{\alpha }^{+}} \sqrt{2(e-(1-\cos \alpha ))} \, \dee \alpha = 
n \hbar 
\end{displaymath}
and $\sqrt{2(e-(1-\cos \alpha ))} > 0$ on $(-{\alpha }^{+}, {\alpha }^{+})$, we obtain 
{\tiny $\left\{ \begin{array}{rl} 
\hspace{-5pt} 0 \le 2m \le n , & \hspace{-8pt}\mbox{if $\varepsilon =1$} \\
\hspace{-5pt} 0 \le 2m -1 \le n , & \hspace{-8pt} \mbox{if $\varepsilon = -1$.} 
\end{array} \right. $}  
\rule{0pt}{15pt}
\indent We now show that 
\begin{equation}
\left\{ \begin{array}{ll} 
n = 2m, & \mbox{if $\varepsilon =1$} \\
n = 2m-1, & \mbox{if $\varepsilon = -1$}. 
\end{array} \right. 
\label{eq-news6six}
\end{equation}
Consider the function 
\begin{displaymath}
F({\alpha }_0) = \frac{1}{2\pi } \int^{{\alpha }_0}_{-{\alpha }_0} \sqrt{2(e-(1-\cos \alpha ))} \, \dee \alpha . 
\end{displaymath}
For ${\alpha }_0 =0$ we get $F(0) =0$; while for ${\alpha }_0 = {\alpha }^{+}$, we get 
$F({\alpha }^{+}) = n \hbar $. By continuity, for every $0 \le k \le n$ there is an angle 
${\alpha }_k \in [0, {\alpha }^{+}]$ such that $F({\alpha }_k) = k \hbar $. Hence 
(\ref{eq-news6six}) holds.  \end{proof} \medskip 

\noindent \textbf{Corollary 6.1.5} \textit{Let $\sigma |_C$ be a nonzero even 
or odd quantum state in 
${\mathfrak{H}}^{\mathrm{even}} \cap {\mathfrak{H}}_0$ or 
${\mathfrak{H}}^{\mathrm{odd}} \cap {\mathfrak{H}}_0$, respectively. Then the 
quantum number of $\sigma |_C$ is the number of ${\Z }_2$-orbits of 
$\varepsilon \, \mu $ $\mathrm{(\ref{eq-nws6one})}$ on the image of the horizontal lift of the curve $\gamma $, which parametrizes $C$.} \medskip    

\begin{proof} Suppose that for some ${\alpha }_0 \in [-{\alpha }^{+}, {\alpha }^{+}]$ we have 
$\widetilde{\gamma }(-{\alpha }_0) = \varepsilon \mu \big( \widetilde{\gamma }({\alpha }_0) \big) $. 
Then $z\big( \gamma (-{\alpha }_0)\big) = \varepsilon z \big( \gamma ({\alpha }_0) \big)$. Repeating the argument of theorem 8 which proves equations (\ref{eq-news6five}) and 
(\ref{eq-news6six}) shows that the quantum number $n$ of the Bohr-Sommerfeld torus 
$C = (H^{\times })^{-1}(e)$ is equal to the number of ${\alpha }_k \in [0, {\alpha }^{+}]$ such 
that $z\big( \gamma ({\alpha }_k) \big) = \varepsilon z\big( \gamma (-{\alpha }_k) \big)$. 
In other words, $n$ is the number of ${\Z }_2$ orbits of $\varepsilon \, \mu $ on the image of 
$\widetilde{\gamma }$. 
\end{proof} \medskip

It follows from lemma 6.1.2 that the operators ${\mathbf{P}}_{\pm }$ enable us to go from 
${\mathfrak{H}}_{+}$ to ${\mathfrak{H}}_{-}$ and back. Thus they play the role 
of shifting operators. In particular ${\{ {\mathbf{P}}_{+}{\sigma }^{+}_m \} }_{m \ge M}$ is a basis 
of ${\mathfrak{H}}_{-}$. Because the parity operator $\mathbf{P}$ induces a ${\Z }_2$-symmetry 
on the Hilbert space $\mathfrak{H}$, by averaging the given inner product, 
we may assume that the parity operator $\mathbf{P}$ preserves the new inner product on 
$\mathfrak{H}$. In order to simplify the presentations we choose the 
orthonormal bases $\{ {\sigma }^{\pm }_m \} $ of ${\mathfrak{H}}_{\pm }$ so that 
${\mathbf{P}}_{\pm}{\sigma }^{\pm }_m = {\sigma }^{\mp}_m$. In order to construct operators 
relating ${\mathfrak{H}}_0$ to ${\mathfrak{H}}_{\pm }$ we need to show that reduction of the 
${\Z}_2$-symmetry of the mathematical pendulum gives rise to 
the quantized ${\Z}_2$-reduced mathematical pendulum. 

\subsection{${\Z }_2$-quantization and reduction}\label{section-6.2} 
%%%%%%%%%%%%%%%%%%%%%%

In this subsection discuss quantization of the ${\Z}_2$-reduced mathematical pendulum. 

\subsubsection{Reduction of the ${\Z }_2$-symmetry}\label{section-6.2.1}
%%%%%%%%%%%%%%%%%%%%%%%%

Here we reduce the ${\Z }_2$-symmetry of the mathematical pendulum 
$(H, T^{\ast }S^1, \omega )$ generated by $\zeta $ (\ref{eq-nwzero}). \medskip 

First we determine the reduced phase space $\widetilde{P}$, which is the 
space $T^{\ast }S^1/{\Z }_2$ of orbits of the ${\Z }_2$-symmetry on $T^{\ast}S^1$. 
To start with we use cut and paste geometric methods to construct the ${\Z }_2$-orbit space. 
Recall that a connected 
subset $\Delta $ of $T^{\ast }S^1$ is a \emph{fundamental domain} for the 
${\Z }_2$-symmetry generated by $\zeta $ (\ref{eq-nwzero}), if it contains 
exactly one point of each ${\Z }_2$-orbit in $T^{\ast }S^1$. \medskip 

\noindent \textbf{Claim 6.2.1} \textit{The set $\Delta = \{ (p, \alpha ) \in T^{\ast }S^1 \, \setrule \, p > 0 \, \, \mathrm{or} \, \, 
p =0 \, \& \, \alpha \in [0, \pi ] \}$ is a fundamental domain for the ${\Z }_2$-action 
generated by $\zeta $.} \medskip 

\begin{proof} Clearly $\Delta $ is connected. Let $(p, \alpha ) \in 
T^{\ast }S^1 \setminus \Delta $. If 
$p \ne 0$, then $p < 0$. So $\zeta (p, \alpha ) = (-p, -\alpha ) \in \Delta $. Suppose that $p =0$ and $\alpha \in (\pi , 2\pi )$. Then $-\alpha \in (0, \pi )$. So $\zeta (0, \alpha ) = 
(0, -\alpha ) \in \Delta $. Hence 
$\zeta (T^{\ast }S^1 \setminus \Delta ) \subseteq \Delta $, which implies 
$T^{\ast }S^1 \setminus \Delta = 
\zeta \big( \zeta \big( T^{\ast }S^1 \setminus \Delta \big) \big) \subseteq \zeta (\Delta )$. Consequently, 
\begin{displaymath}
T^{\ast }S^1 = (T^{\ast }S^1 \setminus \Delta ) \cup \Delta \subseteq \zeta (\Delta ) \cup \Delta  
\subseteq T^{\ast }S^1, 
\end{displaymath}
that is, $\zeta (\Delta ) \cup \Delta  = T^{\ast }S^1$. Note that $\Delta \cap \zeta (\Delta ) 
= \{ (0,0), \, (0, \pi ) \}$, which are the fixed points of the ${\Z}_2$-action. 
\end{proof} \medskip

\noindent Look at the closure $\overline{\Delta}$ of $\Delta $ in $T^{\ast} S^1$ and identify the points on the  boundary of $\overline{\Delta }$, which lie on the same ${\Z}_2$-orbit. The resulting space is a 
model for the orbit space $T^{\ast }S^1/{\Z }_2$. \medskip

We now give another construction for the reduced phase space using invariant theory and 
the concept of a differential space, see \cite{cushman-bates15}. 
The algebra of real analytic functions on $T^{\ast }S^1$, which are invariant under the  
symmetry group ${\Z }_2$, is generated by 
\begin{equation}
{\tau }_1 = \cos \alpha , \qquad {\tau }_2 = p \sin \alpha, \qquad {\tau }_3 = 
\onehalf p^2  - \cos \alpha +1.
\label{eq-s6two}
\end{equation}
These invariant functions are subject to the relation 
\begin{equation}
C(\tau )  = \onehalf {\tau }^2_2 - ({\tau }_3 + {\tau }_1 -1)(1-{\tau }^2_1) =0, 
\quad |{\tau }_1| \le 1 \, \, \& \, \, {\tau }_3 \ge 0,  
\label{eq-s6three}
\end{equation}
which defines the ${\Z }_2$-orbit space $\widetilde{P} = T^{\ast }S^1/{\Z }_2$ as a semialgebraic variety in ${\R }^3$ with coordinates $\tau = ({\tau }_1, {\tau }_2, {\tau }_3)$. We say that a function $f$ on 
$\widetilde{P}$ is \emph{smooth} if there is a smooth function 
$F$ on ${\R }^3$ such that $f = F|_{\widetilde{P}}$. Let \linebreak %
\mbox{}\hspace{1.5in}\begin{tabular}{l} 
\setlength{\unitlength}{2pt} \\
\includegraphics[width=120pt]{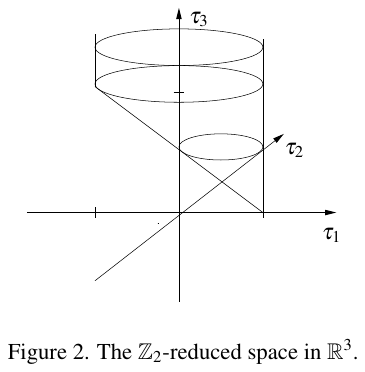}\\
\vspace{-.2in}
\end{tabular}
\par \noindent $C^{\infty}(\widetilde{P})$ be the space of smooth functions on $\widetilde{P}$. Then 
$\big( \widetilde{P}, C^{\infty}(\widetilde{P}) \big) $ is a locally compact subcartesian differential 
space, because $\widetilde{P}$ is a semialgebraic. \medskip 

Next we construct the reduced Hamiltonian. Since the Hamiltonian $H$ of the mathematical 
pendulum is invariant under the ${\Z }_2$-symmetry, it induces a smooth function 
$\widetilde{H}$ on $\widetilde{P}$ given by restricting the smooth function ${\tau }_3: {\R }^3 \rightarrow 
\R :\tau \mapsto {\tau }_3$ to $\widetilde{P}$. \medskip 

In order to have dynamics on $\widetilde{P}$, we first need a Poisson bracket 
${\{ \, \, , \, \, \} }_{{\R }^3}$ on $C^{\infty}({\R }^3)$.  A calculation using the 
Poisson bracket $\{ \, \, , \, \, \} $ on $P$ shows that 
\begin{align}
\{ {\tau }_1, {\tau }_2 \} & =  {\tau}^2_1 -1 = \frac{\partial C}{\partial {\tau }_3} \notag \\
\{ {\tau }_2, {\tau }_3 \} & = 2{\tau }_1({\tau }_3 + {\tau }_1 -1) +{\tau }^2_1 - 1 = \frac{\partial C}{\partial {\tau }_1} \notag \\
\{ {\tau }_3 , {\tau }_1 \} & = {\tau }_2 =   \frac{\partial C}{\partial {\tau }_2} \notag 
\end{align}
For every $F$, $G \in C^{\infty}({\R }^3)$ let 
\begin{equation} 
{\{ F, G \} }_{{\R }^3} = \sum_{i,j} \frac{\partial F}{\partial {\tau }_i} 
\frac{\partial G}{\partial {\tau }_j} \{ {\tau }_i, {\tau }_j \} = \langle \mathrm{grad}\, F \times 
\mathrm{grad}\, G, \mathrm{grad}\, C \rangle ,
\label{eq-news6four}
\end{equation}
Here $\langle \, \, , \, \, \rangle $ is the Euclidean inner product on ${\R }^3$ and $\times $ is the 
vector product. Then ${\{ \, \, , \, \, \} }_{{\R }^3}$ is a Poisson bracket on $C^{\infty}({\R }^3)$. 
On $C^{\infty}(\widetilde{P})$ define a Poisson bracket ${\{ \, , \, \}}_{\widetilde{P}}$ as follows. Suppose that $f, g \in C^{\infty}(\widetilde{P})$. Then there are 
$F,G \in C^{\infty}({\R }^3)$ such that 
$f = F|_{\widetilde{P}}$ and $g = G|_{\widetilde{P}}$. Let ${\{ f, g \}}_{\widetilde{P}} = 
{\{ F, G \}}_{{\R }^3}|_{\widetilde{P}}$. Because of (\ref{eq-news6four}), the defining 
function $C$ (\ref{eq-s6three}) of $\widetilde{P}$ is a Casimir in the Poisson algebra 
$\mathcal{A} = \big( C^{\infty}({\R }^3), {\{ \, , \, \}}_{{\R }^3}, \cdot \big)$. 
Hence, the collection $\mathcal{I}$ of all smooth functions on ${\R }^3$, which vanish identically on $\widetilde{P}$, is a Poisson ideal in $\mathcal{A}$. Consequently, the Poisson bracket 
${\{ \, , \, \}}_{\widetilde{P}}$ is well defined and $\mathcal{B} = \mathcal{A}/\mathcal{I} = 
\big( C^{\infty}(\widetilde{P}) = 
C^{\infty}({\R }^3)|\mathcal{I},{\{ \, , \, \}}_{\widetilde{P}}, \cdot \big)$ is a Poisson algebra. \medskip 

Consider the derivation $-{\ad }_{{\tau }_3}$ on the Poisson algebra 
$\mathcal{A}$. This derivation gives rise to the ${\Z }_2$-reduced Hamiltonian vector field 
$-{\ad}_{\widetilde{H}}$ on the locally compact subcartesian differential space 
$(\widetilde{P}, C^{\infty}(\widetilde{P}))$ associated to the ${\Z }_2$-reduced Hamiltonian 
$\widetilde{H}$. To see this note that on ${\R }^3$ the integral curves of $-{\ad }_{{\tau }_3}$ satisfy
\begin{align}
{\dot{\tau }}_1 & =  {\{ {\tau }_1, {\tau }_3 \}}_{{\R }^3}  = - {\tau }_2 \notag \\
{\dot{\tau }}_2 & =  {\{ {\tau }_2, {\tau }_3 \} }_{{\R }^3}  = 
2{\tau }_1({\tau }_3 + {\tau }_1 -1) +{\tau }^2_1-1 \notag \\
{\dot{\tau }}_3 & = {\{ {\tau }_3, {\tau }_3 \} }_{{\R }^3}  = 0. \notag 
\end{align}
Because $C$ is a Casimir of the Poisson algebra $\mathcal{A}$, we obtain 
$0 = {\{ C, {\tau }_3 \} }_{{\R }^3}$. In other words, $C$ is an integral of $-{\ad }_{{\tau }_3}$. 
A calculation shows that $-{\ad }_{{\tau }_3}$ leaves the sets $C^{-1}(0)$, 
$\{ {\tau }_3 +{\tau }_1 -1=0 \} $, and $\{ {\tau }_1 = \pm 1 \}$ invariant. Thus the ${\Z }_2$-reduced space 
$\widetilde{P}$ is invariant under the flow of $-{\ad }_{{\tau }_3}$. Consequently, the ${\Z }_2$-reduced Hamiltonian vector field $-{\ad }_{\widetilde{H}}$, where $\widetilde{H} = {\tau }_3|_{\widetilde{P}}$, is 
defined on $\widetilde{P}$. Because the Hamiltonian vector field $X_H$ of the mathematical pendulum is complete, the reduced vector field $-{\ad }_{\widetilde{H}}$ is complete. Its flow 
${\varphi }^{\widetilde{H}}_t$ is a $1$-parameter group of diffeomorphisms of 
$\widetilde{P}$. In fact, for $\widetilde{p} \in {\widetilde{H}}^{-1}(e)$ the closure of the integral curve 
$t \mapsto {\varphi }^{\widetilde{H}}_t(\widetilde{p})$ is a connected component of the level set 
${\widetilde{H}}^{-1}(e)$, since a level set of the reduced Hamiltonian $\widetilde{H}$ is compact.\medskip  

We now take a closer look at the ${\Z }_2$-reduction mapping 
\begin{equation}
\widetilde{\pi }: T^{\ast }S^1 \rightarrow \widetilde{P} \subseteq {\R }^3: (p, \alpha ) \mapsto 
\tau (p, \alpha ) .
\label{eq-nws6fourstarstar}
\end{equation}
The ${\Z }_2$-action on $T^{\ast }S^1$ has two fixed points: $p_0 = (0,0)$ and 
$p_2 = (0, \pi )$. So the reduced space $\widetilde{P}$ has two singular points 
${\widetilde{p}}_0 = (1,0,0)$ and ${\widetilde{p}}_2 = (-1,0,2)$, which are conical. 
The set ${\widetilde{P}}^{\times }$ of nonsingular points of $\widetilde{P}$ is 
$\widetilde{P}\setminus \{ {\widetilde{p}}_0, {\widetilde{p}}_2 \} $, which is a smooth 
manifold that is diffeomorphic to ${\R }^2 \setminus \{ (\pm 1, 0 ) \}$. The ${\Z}_2$-orbit 
map $\widetilde{\pi }$ (\ref{eq-nws6fourstarstar}) restricted to $T^{\times }S^1 = T^{\ast }S^1 \setminus 
\{ p_0, p_2 \} $ is the proper submersion 
\begin{equation}
{\widetilde{\pi }}^{\times }: T^{\times }S^1 \rightarrow {\widetilde{P}}^{\times }: 
(p, \alpha ) \mapsto \tau (p, \alpha ), 
\label{eq-nws6doublestar}
\end{equation}
whose fiber $({\widetilde{\pi}}^{\times})^{-1}(\tau )$ at $\tau \in \widetilde{P}$ is two distinct points. Thus we have proved \medskip 

\noindent \textbf{Lemma 6.2.2} \textit{The ${\Z}_2$-orbit map ${\widetilde{\pi }}^{\times }$ 
$\mathrm{(\ref{eq-nws6doublestar})}$ is a $2$ to $1$ covering map.} \medskip   

The $1$-form ${\theta }^{\times } = \theta |_{T^{\times }S^1} = (p \, \dee \alpha )|_{T^{\times }S^1}$ on $T^{\times} S^1$ is invariant under the ${\Z }_2$-action, since $T^{\times }S^1$ is a smooth 
${\Z}_2$-invariant manifold. Hence ${\theta }^{\times }$ pushes down under the 
${\Z }_2$-orbit map ${\widetilde{\pi }}^{\times }$ (\ref{eq-nws6doublestar}) to a 
$1$-form ${\widetilde{\theta }}^{\times }$ on the smooth manifold 
${\widetilde{P}}^{\times }$. 
So $({\widetilde{\pi }}^{\times })^{\ast }{\widetilde{\theta }}^{\times } = {\theta }^{\times }$. 
Here are explicit expressions for the $1$-form ${\widetilde{\theta }}^{\times }$. \medskip   

\noindent \textbf{Claim 6.2.3} \textit{On ${\widetilde{U}}_1 = {\widetilde{P}}^{\times } \setminus 
\{ {\tau }_1 = \pm 1 \}$ we have ${\widetilde{\theta }}^{\times }|{\widetilde{U}}_1 = 
- {\tau }_2(1-{\tau }^2_1)^{-1}\, \dee {\tau }_1$; while on 
${\widetilde{U}}_2 = {\widetilde{P}}^{\times } \setminus ( \{ {\tau }_1 =0 \} 
\cup \{ {\tau }_3 +{\tau }_1 -1 =0 \} )$ we have ${\widetilde{\theta }}^{\times } |{\widetilde{U}}_2 = 
\big( 2({\tau }_3 +{\tau }_1 -1) \dee {\tau }_2 -{\tau }_2 \dee {\tau }_1 - 
{\tau }_2 \dee {\tau }_3 \big) 
\big( 2{\tau }_1({\tau }_3 +{\tau }_1 -1) \big)^{-1}$. Note ${\widetilde{P}}^{\times } = 
{\widetilde{U}}_1 \cup {\widetilde{U}}_2$.} \vspace{-.1in}  

\begin{proof} On $U_1 = T^{\times }S^1 \setminus \{ (p, \pm \pi ) \in T^{\times }S^1 \, \setrule \, p \in {\R }^{\times } \}$ we have 
\begin{align}
({\widetilde{\pi }}^{\times })^{\ast }( {\widetilde{\theta }}^{\times }|{\widetilde{U}}_1) & = 
-\frac{p \sin \alpha }{1 - {\cos }^2 \alpha} \dee (\cos \alpha ) = 
(p \, \dee \alpha )|U_1 = {\theta }^{\times }|U_1 ; \notag 
\end{align} 
while on $U_2 = T^{\times }S^1 \setminus \{ (p, \pm \pi /2 ) \in T^{\times }S^1 \, \setrule \, p \in 
{\R }^{\times }\} \cup \{ (0, \alpha  ) \in T^{\times }S^1 \} $ we have 
\begin{align} 
({\widetilde{\pi }}^{\times })^{\ast } ({\widetilde{\theta }}^{\times } |{\widetilde{U}}_2) & = 
\frac{\big( p^2 \dee (p \sin \alpha ) -p\sin \alpha \, \dee (\cos \alpha ) 
- p\sin \alpha (p\, \dee p - \dee (\cos \alpha )) \big) }{p^2 \cos \alpha } \notag \\
& = (p \, \dee \alpha )|U_2 = {\theta }^{\times }|U_2 . \notag 
\end{align}
Note that for $i =1$, $2$ we have $({\widetilde{\pi }}^{\times })^{-1}({\widetilde{U}}_i) = U_i$ and $T^{\times }S^1 = U_1 \cup U_2$. 
\end{proof} \medskip 

\noindent Because the $2$-form $\omega = \dee \theta $ on $T^{\ast }S^1$ is invariant 
under the ${\Z }_2$-symmetry generated by $\zeta $ (\ref{eq-nwzero}), the $2$-form 
${\omega }^{\times } = {\omega }|_{T^{\times }S^1} = \dee (\theta |_{T^{\times }S^1} ) = 
\dee {\theta }^{\times }$ on $T^{\times}S^1$ is ${\Z }_2$-invariant. Hence ${\omega }^{\times }$ pushes down to a $2$-form ${\widetilde{\omega }}^{\times }$ on the ${\Z }_2$-reduced space 
${\widetilde{P}}^{\times }$. Now 
\begin{displaymath}
({\widetilde{\pi }}^{\times })^{\ast }(\dee {\theta }^{\times }) = 
\dee ( ({\widetilde{\pi }}^{\times })^{\ast }{\widetilde{\theta }}^{\times }) = 
\dee {\theta }^{\times } = {\omega }^{\times } = 
({\widetilde{\pi }}^{\times })^{\ast } {\widetilde{\omega }}^{\times }.
\end{displaymath}
Since ${\widetilde{\pi }}^{\ast }$ is surjective, we obtain ${\widetilde{\omega }}^{\times } = 
\dee {\widetilde{\theta }}^{\times }$. Thus the punctured ${\Z }_2$-reduced space 
${\widetilde{P}}^{\times }$ is a symplectic manifold with symplectic form 
${\widetilde{\omega }}^{\times } = \dee {\widetilde{\theta }}^{\times }$. \medskip 

We now compute the ${\Z }_2$-reduced actions. On ${\widetilde{P}}^{\times }
\setminus \{ {\tau }_1 = \pm 1 \}$ the 
$1$-form ${\widetilde{\theta }}^{\times }$ is $ -{\tau }_2(1-{\tau }^2_1)^{-1} \, \dee {\tau }_1$, where ${\tau }_2 = \mp \sqrt{2({\tau }_3 +
{\tau }_1 -1)(1-{\tau }^2_1)}$. So  ${\widetilde{\theta }}^{\times } = \pm \sqrt{\frac{2({\tau }_3 +
{\tau }_1 -1)}{1-{\tau }^2_1}}\, \dee {\tau }_1$. The ${\Z }_2$-reduced Hamiltonian on 
${\widetilde{P}}^{\times }$ is ${\widetilde{H}}^{\times } = \widetilde{H}|{\widetilde{P}}^{\times }$. Consequently, the reduced action ${\widetilde{I}}^{\times }: (0,2) \cup (2, \infty ) \rightarrow \R $ on a connected component $\widetilde{C}(e)$ of the level set 
$({\widetilde{H}}^{\times })^{-1}(e)$ is 
\begin{align}
{\widetilde{I}}^{\times }(e)& = \frac{1}{2\pi} \, \int_{\widetilde{C}(e)} {\widetilde{\theta }}^{\times } = 
\frac{1}{\pi} \int^1_{\max (1-e, -1)} \frac{\sqrt{2(e+{\tau }_1 -1)}}{\sqrt{1-{\tau }^2_1}} \dee {\tau }_1,
\label{eq-onenew*}
\end{align} 
when $0 < e < 2$ or $e >2$. We now calculate the integral in (\ref{eq-onenew*}). First we consider the case when $0 < e <2$. Letting $u^2 = {\tau }_1 -(1-e)$,  
$u = \sqrt{e}v$, and then $v = \cos \varphi $, we get successively
\begin{align}
{\widetilde{I}}^{\times }(e) & =  
\frac{\sqrt{2}}{\pi} \int^{\sqrt{e}}_0 \frac{2u^2}{\sqrt{(e-u^2)(2-e+u^2)}} \, \dee u \notag  \\ 
& = \frac{2\sqrt{2}}{\pi} e \, \int^1_0 \frac{v^2}{\sqrt{(1-v^2)(2-e +ev^2)}} \, \dee v \notag \\
& = \frac{2}{\pi }e \int^{\pi /2}_0 \frac{{\cos }^2 \varphi }{\sqrt{1-\frac{e}{2}{\sin }^2{\varphi }}} \, 
\dee \varphi  \ge 0.
\label{eq-twonew*}
\end{align}
Next we treat the case when $e > 2$. Letting ${\tau }_1 = \cos \vartheta $ and $\vartheta  = 2 \varphi $ successively, 
we get 
\begin{align}%
{\widetilde{I}}^{\times }(e) &  = \frac{\sqrt{2}}{\pi } \int^{\pi /2}_0 \sqrt{e-1+\cos \vartheta } \, \dee \vartheta 
= \frac{2\sqrt{2}}{\pi } \, \int^{\pi /2}_0 \sqrt{e-2 +2\, {\cos }^2\varphi } \, \dee \varphi  \notag \\
& = \frac{2\sqrt{2e}}{\pi } \, \int^{\pi /2}_0 \sqrt{1 - \ttfrac{2}{e} \, {\sin }^2\varphi } \, \dee \varphi \ge 0. 
\label{eq-threenew*}
\end{align}
Note that the ${\Z }_2$-reduced action ${\widetilde{I}}^{\times }$ (\ref{eq-twonew*}) and 
(\ref{eq-threenew*}) is one half the original action $I^{\times} = I|_{T^{\times }S^1}$ 
(\ref{eq-s1threedot}) and (\ref{eq-s1threedagg}). Moreover, the function 
$e \mapsto {\widetilde{I}}^{\times }(e)$ is continuous at $e = 2$. 
Dullin \cite{dullin} shows, ${\widetilde{I}}^{\times }$ has a logarithm term in its series expansion 
in $e-2$, which shows ${\widetilde{I}}^{\times }$ is not differentiable at $e =2$. \medskip 

The corresponding ${\Z }_2$-reduced angle ${\widetilde{\vartheta}}^{\times }$ is 
\begin{equation}
{\widetilde{\vartheta }}^{\times } = \frac{2\pi }{\widetilde{T}}t = 
\frac{2\pi }{\widetilde{T}} \int^1_{t} \frac{2 \, \dee {\tau }_1}{\sqrt{2(1-{\tau }^2_1)(e-1+{\tau }_1)}}, 
\label{eq-threestarnew}
\end{equation}
where $t \in [ \max (1-e, -1), 1]$ and 
\begin{equation}
\widetilde{T} = \widetilde{T}(e) = \int^1_{-1} \frac{2 \, \dee {\tau }_1}{\sqrt{2(1-{\tau }^2_1)(e-1+{\tau }_1)}}.
\label{eq-nws6threedoublestar}
\end{equation} 

In order to simplify the notation, in the following we will use $\widetilde{I}$ for the reduced action 
on the reduced space $\widetilde{P}$ and $\widetilde{\vartheta}$ for the reduced angle, which is not defined at the singular points ${\widetilde{p}}_0$ and ${\widetilde{p}}_2$ of $\widetilde{P}$. 
Note that ${\widetilde{I}}^{\times} = \widetilde{I}|_{{\widetilde{P}}^{\times }}$ and 
${\widetilde{\vartheta }}^{\times } = \widetilde{\vartheta}|_{{\widetilde{P}}^{\times }}$. 

\subsubsection{Reduction of the ${\Z }_2$-quantum symmetry}\label{section-6.2.2}
%%%%%%%%%%%%%%%%%%%%%%%%

In this subsubsection we reduce the ${\Z }_2$-symmetry of the quantized mathematical 
pendulum. In other words, we reduce the ${\Z}_2$-action on $L \subseteq \C \times {\R }^2$ 
generated by the mapping $\mu $ (\ref{eq-nws6one}). We use invariant theory. \medskip 

The algebra of invariant real analytic functions is generated by 
\begin{displaymath}
{\sigma }_1 = z, \, \, {\tau }_1 = \cos \alpha , \, \, {\tau }_2 = p \sin \alpha , \, \, {\tau }_3 = \onehalf p^2 - \cos \alpha + 1 
\end{displaymath}
subject to the relation  
\begin{align}
\onehalf {\tau }^2_2 & = \onehalf p^2 {\sin }^2 \alpha = \onehalf p^2(1-{\cos }^2 \alpha ) \notag \\
& = ({\tau }_1 +{\tau }_3 -1)(1-{\tau }^2_1), \quad |{\tau }_1| \le 1 \, \, \& \, \, {\tau }_3 \ge 0.  
\label{eq-nwseven} 
\end{align}
Equation (\ref{eq-nwseven}) defines the ${\Z}_2$-orbit space 
$\widetilde{\mathcal{P}} = L/{\Z }_2$. The Hilbert mapping 
\begin{displaymath}
\widetilde{\varsigma }: L \rightarrow 
\widetilde{\mathcal{P}}: \big( z, (p, \alpha ) \big) \mapsto \big( {\sigma }_1(z), {\tau }_1(p,\alpha ), 
{\tau }_2(p, \alpha ), {\tau }_3(p,\alpha ) \big) = \big( {\sigma }_1(z), \tau (p,\alpha ) \big)
\end{displaymath} 
is the orbit map of the ${\Z }_2$-action. The ${\Z}_2$-orbit space 
$\widetilde{\mathcal{P}}$ is $\C \times (T^{\ast }S^1/{\Z }_2)$, which is a semialgebraic 
variety with two singular planes $\C \times \{ (\pm 1, 0, 1\mp 1) \} $. We view 
$\widetilde{\mathcal{P}}$ as a complex ``line bundle'' over the ${\Z}_2$-orbit space 
$\widetilde{P}$ with bundle projection
\begin{displaymath}
\widetilde{\varpi }:\widetilde{\mathcal{P}} = \C \times \widetilde{P} \rightarrow \widetilde{P}:
 \big( {\sigma }_1, \tau = ({\tau }_1, {\tau }_2, {\tau }_3) \big) \mapsto \tau .
\end{displaymath}
Here (\ref{eq-nwseven}) is the defining relation of the ${\Z }_2$-orbit space $\widetilde{P}$ with orbit mapping $\widetilde{\pi }: T^{\ast }S^1 \rightarrow \widetilde{P}: (p, \alpha ) \mapsto \tau ( p, \alpha )$. 
Consider the smooth manifold ${\widetilde{\mathcal{P}}}^{\times} = 
\widetilde{\mathcal{P}} \setminus \{ (\pm 1, 0, 1\mp 1) \} = \C \times {\widetilde{P}}^{\times }$ of nonsingular points of $\widetilde{\mathcal{P}}$, where  
${\widetilde{P}}^{\times} = T^{\times }S^1/{\Z}_2$ and $T^{\times }S^1 = T^{\ast }S^1 \setminus 
\{ (0, 0), (0, \pi ) \}$. Restricting $\widetilde{\varpi}$ to $\C \times {\widetilde{P}}^{\times }$ gives a trivial smooth bundle ${\widetilde{\rho}}^{\, \times }: {\widetilde{L}}^{\times} = \C \times {\widetilde{P}}^{\times} \rightarrow {\widetilde{P}}^{\times}$. The bundle ${\widetilde{\rho }}^{\, \times }$ serves as the 
quantum bundle of the ${\Z}_2$-reduced mathematical pendulum 
$({\widetilde{H}}^{\times }, {\widetilde{P}}^{\times }, {\widetilde{\omega }}^{\times })$. 

\subsubsection{The ${\Z }_2$-reduced quantum system}\label{section-6.2.3}
%%%%%%%%%%%%%%%%%%%%%

In this subsubsection we quantize the ${\Z }_2$-reduced quantum mathematical pendulum, 
namely, the ${\Z }_2$-reduced mathematical pendulum $({\widetilde{H}}^{\times }, 
{\widetilde{P}}^{\times } , {\widetilde{\omega }}^{\times } )$ with the ${\Z }_2$-reduced 
quantum bundle 
\begin{displaymath}
{\widetilde{\rho }}^{\, \times }: {\widetilde{L}}^{\times } = \C \times {\widetilde{P}}^{\times } 
\rightarrow {\widetilde{P}}^{\times }: (z, \tau ) \mapsto \tau 
\end{displaymath}
and trivializing section ${\widetilde{\lambda }}_0: {\widetilde{P}}^{\times } \rightarrow 
{\widetilde{L}}^{\times }: \tau \mapsto (1, \tau )$. \medskip 

We need to find a connection ${\widetilde{\nabla }}^{\times }$ on the smooth sections of the bundle ${\widetilde{\rho }}^{\, \times }$, which is related to the original connection 
${\nabla }^{\times }$ on $P^{\times }$. On smooth sections of the bundle ${\rho }^{\times }$ 
we have a connection, whose covariant derivative ${\nabla }^{\times }_X$ in the direction 
of the smooth vector field $X$ on $T^{\times }S^1$ acts on the section ${\lambda }^{\times }_0$ by ${\nabla }^{\times }_{X}{\lambda }^{\times }_0 = 
-i{\hbar}^{-1} \big( X \lefthook {\theta }^{\times } \big) 
{\lambda }^{\times}_0$. Suppose that $\widetilde{X}$ is a smooth vector field on 
${\widetilde{P}}^{\times }$, which is ${\widetilde{\pi }}^{\times }$-related to the vector field $X$, that is, $T{\widetilde{\pi }}^{\times } X  = \widetilde{X} \comp {\widetilde{\pi }}^{\times }$. On the line bundle ${\widetilde{\rho }}^{\, \times }$ with trivializing section 
${\widetilde{\lambda}}^{\times }_0$ 
define a connection ${\widetilde{\nabla }}^{\times }$ by  $({\widetilde{\pi }}^{\times})^{\ast }
({\widetilde{\nabla }}^{\times }_{\widetilde{X}} {\widetilde{\lambda }}^{\times }_0) = 
{\nabla }^{\times }_X {\lambda }^{\times }_0$. In other words, \medskip

\noindent \textbf{Fact 6.2.3.1}
\begin{equation}
{\widetilde{\nabla }}^{\times }_{\widetilde{X}} {\widetilde{\lambda }}^{\times }_0 = 
-i{\hbar }^{-1} (\widetilde{X} \lefthook {\widetilde{\theta }}^{\times }) 
{\widetilde{\lambda }}^{\times }_0.
\label{eq-s3newone}
\end{equation}

\begin{proof} Equation (\ref{eq-s3newone}) follows because by definition
\begin{displaymath}
({\widetilde{\pi }}^{\times })^{\ast }\big( {\widetilde{\nabla }}^{\times }_{\widetilde{X}} 
{\widetilde{\lambda }}^{\times }_0 \big) = 
{\nabla }^{\times }_X {\lambda }^{\times }_0 = -i{\hbar }^{-1}(X \lefthook {\theta }^{\times }) {\lambda }^{\times}_0; 
\end{displaymath}
whereas 
\begin{align}
({\widetilde{\pi }}^{\times})^{\ast }\big( -i{\hbar}^{-1}(\widetilde{X} \lefthook {\widetilde{\theta }}^{\times })
{\widetilde{\lambda }}^{\times }_0 \big)& = 
-i{\hbar}^{-1}\big( (\widetilde{X} \lefthook {\widetilde{\theta }}^{\times }) \comp 
{\widetilde{\pi }}^{\times} \big) 
{\widetilde{\lambda }}^{\times }_0 \comp {\widetilde{\pi }}^{\times} \notag \\
&\hspace{-1in} = -i{\hbar}^{-1}(\widetilde{X} \comp {\widetilde{\pi }}^{\times} \lefthook 
{\widetilde{\theta }}^{\times } 
\comp {\widetilde{\pi }}^{\times} ) {\widetilde{\lambda }}^{\times }_0 \comp 
{\widetilde{\pi }}^{\times} 
\notag \\
&\hspace{-1in} = -i{\hbar}^{-1}(T{\widetilde{\pi }}^{\times} X \lefthook 
T{\widetilde{\pi }}^{\times} {\theta }^{\times }) 
{\widetilde{\lambda }}^{\times }_0 \comp {\widetilde{\pi }}^{\times} \notag \\
& \hspace{-1in} = -i{\hbar}^{-1}( X \lefthook {\theta }^{\times}) \big( (T{\widetilde{\pi }}^{\times})^t {\widetilde{\lambda }}^{\times }_0 \comp {\widetilde{\pi }}^{\times} \big) \notag \\
& \hspace{-1in}= -i{\hbar}^{-1}(X \lefthook {\theta }^{\times }) {\lambda }^{\times }_0, \quad 
\mbox{since $({\widetilde{\pi }}^{\times})^{\ast }{\widetilde{\lambda }}^{\times }_0 = 
{\lambda }^{\times }_0$.}  
\notag 
\end{align}
Thus $({\widetilde{\pi }}^{\times})^{\ast } \big( {\widetilde{\nabla }}^{\times }_{\widetilde{X}} {\widetilde{\lambda }}_0 \big) = 
({\widetilde{\pi }}^{\times})^{\ast }\big( -i{\hbar}^{-1}(\widetilde{X} \lefthook 
{\widetilde{\theta }}^{\times }){\widetilde{\lambda }}_0 \big) $, 
which implies (\ref{eq-s3newone}) since ${\widetilde{\pi }}^{\times}$ is surjective. 
\end{proof} \medskip

Next we determine the quantization rules for the ${\Z}_2$-reduced Hamiltonian system 
$( {\widetilde{H}}^{\times } = \widetilde{H}|{\widetilde{P}}^{\times }, {\widetilde{P}}^{\times}, {\widetilde{\omega }}^{\times} )$ with quantum line bundle ${\widetilde{\rho }}^{\, \times }$ and 
trivializing section ${\widetilde{\lambda}}_0$. 
The mapping ${\widetilde{P}}^{\times } \rightarrow T{\widetilde{P}}^{\times }: p \mapsto \spann \{ X_{{\widetilde{H}}^{\times }}(p) \} $ defines a smooth Lagrangian distribution $\widetilde{D}$ on the symplectic manifold $({\widetilde{P}}^{\times }, {\widetilde{\omega }}^{\times })$, which is a polarization of $({\widetilde{P}}^{\times },{\widetilde {\omega }}^{\times})$. A leaf of $\widetilde{D}$ is a connected component of a level set $({\widetilde{H}}^{\times })^{-1}(e)$ of the ${\Z }_2$-reduced Hamiltonian 
${\widetilde{H}}^{\times }$ on ${\widetilde{P}}^{\times }$, which is a smooth $S^1$ when 
$e \in (0, 2) \cup (2, \infty )$. \medskip  

Let $\gamma : \R \rightarrow {\widetilde{P}}^{\times }$ be an integral curve of 
$X_{{\widetilde{H}}^{\times }}$ of energy $e \in (0,2) \cup (2, \infty )$. 
Then $\gamma $ is periodic of primitive period $\widetilde{T} = \widetilde{T}(e) >0$. Also $\gamma $ parametrizes a connected component ${\widetilde{C}}^{\times }(e)$ of the 
smooth level set $({\widetilde{H}}^{\times })^{-1}(e)$. Parallel transport the section 
$\widetilde{\lambda }$ of the $C^{\times}$-bundle ${\widetilde{\, \rho }}^{\times }$ along 
$\gamma $ using the connection ${\widetilde{\nabla }}^{\times }$. 
Then at every point $\gamma (t)$ in $({\widetilde{H}}^{\times })^{-1}(e)$ we have 
\begin{displaymath}
0 = ({{\widetilde{\nabla }}^{\times }}_{X_{{\widetilde{H}}^{\times }}} \widetilde{\lambda } ) (\gamma (t)) = (L_{X_{{\widetilde{H}}^{\times }}}f)(\gamma (t)){\widetilde{\lambda }}_0 - 
\ttfrac{i}{\hbar} (X_{{\widetilde{H}}^{\times }} \lefthook 
{\widetilde{\theta }}^{\times })(\gamma (t)) f(\gamma (t)){\widetilde{\lambda }}_0, 
\end{displaymath}
that is, 
\begin{equation}
\frac{\dee F(t)}{\dee t} - \frac{i}{\hbar } (X_{{\widetilde{H}}^{\times }} \lefthook 
{\widetilde{\theta }}^{\times })
(\gamma (t)) F(t) = 0 , 
\label{eq-s2one}
\end{equation}
where $F(t) = f(\gamma (t))$. For equation (\ref{eq-s2one}) to have a nonvanishing solution 
\begin{displaymath}
F(\widetilde{T}) = F(0) \exp \Big( \ttfrac{i}{\hbar } \int^{\widetilde{T}}_0 
(X_{{\widetilde{H}}^{\times}} \lefthook {\widetilde{\theta }}^{\times } )
(\gamma (t)) \, \dee t \Big) = A(\widetilde{T}) F(0), 
\end{displaymath}
the holonomy $A(\widetilde{T})$ of the connection ${\widetilde{\nabla }}^{\times }$ along $\gamma $ 
must equal $1$, because $F(\widetilde{T}) =f(\gamma (\widetilde{T})) = f(\gamma (0)) 
= F(0) \ne 0$. Consequently, for some $k \in \Z $ we have 
\begin{align}
k & = \frac{1}{\hbar } \int^{\widetilde{T}}_0 (X_{{\widetilde{H}}^{\times }} \lefthook 
{\widetilde{\theta }}^{\times } )(\gamma (t)) \, \dee t  = 
\frac{1}{\hbar } \int^{\widetilde{T}}_0 {\widetilde{\theta }}^{\times } (\gamma (t))
X_{{\widetilde{H}}^{\times }}(\gamma (t))  \, \dee t \notag \\
& = \frac{1}{\hbar } \int^{\widetilde{T}}_0 {\widetilde{\theta }}^{\times } (\gamma (t))
\frac{\dee \gamma (t)}{\dee t}  \, \dee t 
= \frac{1}{\hbar } \int^{\widetilde{T}}_0 {\gamma }^{\ast } {\widetilde{\theta }}^{\times }  
= \frac{1}{\hbar } \int_{({\widetilde{H}}^{\times })^{-1}(e)} \hspace{-.25in} {\widetilde{\theta }}^{\times } \notag 
\end{align}
In other words, when $e \in (0,2) \cup (2, \infty) $ the quantization rule for the 
$({\Z }_2, \cdot )$-reduced quantized Hamiltonian system $({\widetilde{H}}^{\times}, 
{\widetilde{P}}^{\times}, {\widetilde{\omega }}^{\times})$ with quantum bundle 
${\widetilde{\rho}}^{\, \times }: {\widetilde{L}}^{\times } \rightarrow {\widetilde{P}}^{\times }$ is 
\begin{equation}
0 \le {\widetilde{I}}^{\times }(e) = 
\frac{1}{2\pi }\int_{{\widetilde{C}}^{\times }(e)} {\widetilde{\theta }}^{\times }  = 
k \hbar , \quad \mbox{for some $k \in {\Z }_{\ge 0}$,}
\label{eq-s2two}
\end{equation} 
where ${\widetilde{I}}^{\times }$ is the action (\ref{eq-onenew*}) of the ${\Z }_2$-reduced mathematical 
pendulum. \medskip 

\noindent \textbf{Lemma 6.2.3.2} \textit{The reduction mapping 
${\widetilde{\pi }}^{\times }:T^{\times }S^1 
\rightarrow {\widetilde{P}}^{\times }: (p, \alpha ) \mapsto \tau (p, \alpha )$ 
$\mathrm{(\ref{eq-nws6doublestar})}$ maps a Bohr-Sommerfeld torus $C(e)$ of the mathematical pendulum $(H^{\times }, T^{\times }S^1, {\omega }^{\times })$ onto a Bohr-Sommerfeld torus $\widetilde{C}(e)$ of the ${\Z }_2$-reduced mathematical pendulum $({\widetilde{H}}^{\times }, {\widetilde{P}}^{\times }, {\widetilde{\omega }}^{\times })$.} \medskip 

\begin{proof} By lemma 6.2.3.1 the ${\Z}_2$-reduction mapping ${\widetilde{\pi }}^{\times }$ is 
a $2$ to $1$ covering map. Its preimage of the $e$-level set $({\widetilde{H}}^{\times })^{-1}(e)$ 
of the ${\Z }_2$-reduced Hamiltonian on the ${\Z }_2$-reduced phase space 
${\widetilde{P}}^{\times }$ is $(H^{\times })^{-1}(e)$ if $0 < e <1$ or $e >2$. Thus the 
image of a connected component $C(e)$ of $(H^{\times })^{-1}(e)$ under 
${\widetilde{\pi }}^{\times }$ is $({\widetilde{H}}^{\times })^{-1}(e)$. Since 
\begin{align}
({\widetilde{\pi }}^{\times })^{\ast }\Big( {\widetilde{I}}^{\times }|_{\widetilde{C}(e)} \Big) & = 
({\widetilde{\pi }}^{\times })^{\ast } \Big( \frac{1}{2\pi } \int_{({\widetilde{H}}^{\times })^{-1}(e)} {\widetilde{\vartheta }}^{\times } \Big) \notag \\
&\hspace{-.5in} = \left\{ \begin{array}{rl}
\frac{2}{2\pi} \int_{(H^{\times })^{-1}(e)} {\vartheta }^{\times } , & \mbox{if $0 < e <2$} \\
\rowspace \frac{2}{2\pi} \int_{C_{\pm}(e)} {\vartheta }^{\times } , & \mbox{if $2 < e $} \end{array} \right. 
= \frac{1}{2\pi } \int_{C(e)} {\vartheta }^{\times } = (I^{\times })|_{C(e)}, \notag 
\end{align}
the image under ${\widetilde{\pi }}^{\times }$ of a Bohr-Sommerfeld torus of the mathematical pendulum is a Bohr-Sommerfeld torus of the ${\Z }_2$-reduced mathematical pendulum.\end{proof} \medskip  

For every positive integer $k$, let ${\widetilde{\sigma }}_k$ be a section of the line bundle 
${\widetilde{\rho }}^{\, \times }$, which is supported and covariantly constant on the level set 
$({\widetilde{I}}^{\times })^{-1}( k\hbar )$. As before we add the quantum number $0$, 
which corresponds to a section supported on the singular Bohr-Sommerfeld torus corresponding 
to the singular point $(1,0,0)$ of $\widetilde{P}$. The collection 
${\{ {\widetilde{\sigma }}_k \} }_{k \in {\Z}_{\ge 0}}$ is an orthonormal basis of the space 
$\widetilde{\mathfrak{H}}$ of quantum states of the ${\Z }_2$-reduced mathematical pendulum. \medskip 

Since the quantum states ${\{ {\widetilde{\sigma }}_k \} }_{k \in {\Z}_{\ge 0}}$ are ordered by 
increasing $k$, there exist shifting operators $\widetilde{\mathbf{a}}$ and 
$\widetilde{\mathbf{b}}$
such that 
\begin{equation}
\begin{array}{rl}
\widetilde{\mathbf{b}}\, {\widetilde{\sigma }}_k = {\widetilde{\sigma }}_{k+1}, & 
\mbox{for $k \ge 0$} \\
\rowspace \widetilde{\mathbf{a}}\, {\widetilde{\sigma }}_k = {\widetilde{\sigma }}_{k-1}, & 
\mbox{for $k > 0$ and $\widetilde{\mathbf{a}}\, {\widetilde{\sigma }}_0 =0$. } 
\end{array}
\label{eq-threenew}
\end{equation}
Because the local lattice structure of the set of Bohr-Sommerfeld tori on $\widetilde{P}$ is linear,
the shifting operators $\widetilde{\mathbf{a}}$ and $\widetilde{\mathbf{b}}$ are also well defined 
across the singularitiy at the reduced energy value $e =2$. As before, the operators 
$\widetilde{\mathbf{a}}$ and $\widetilde{\mathbf{b}}$ satisfy the same commutation relations 
as the quantum operators ${\mathbf{Q}}_{{\mathrm{e}}^{-i\widetilde{\vartheta }}}$ and ${\mathbf{Q}}_{{\mathrm{e}}^{i\widetilde{\vartheta }}}$, respectively. 

\subsubsection{Lifting the shifting operators}\label{section-6.2.4}
%%%%%%%%%%%%%%%%%%%

In this subsubsection we use the isomorphism $\mathbf{R}: {\mathfrak{H}}^{\mathrm{even}} 
\rightarrow \widetilde{\mathfrak{H}}$ to lift the shifting operators on $\widetilde{\mathfrak{H}}$ to 
shifting operators on ${\mathfrak{H}}^{\mathrm{even}}$. \medskip 

We define $\mathbf{R}$ as the operator which sends the basis 
$\{ {\sigma }^0_{2k},\, k=1, \ldots , K;\, $ ${\sigma }^{+}_m + {\sigma }^{-}_m, \, m = M, M+1, \ldots   \} $ of ${\mathfrak{H}}^{\mathrm{even}}$ 
to the basis $\{ {\widetilde{\sigma }}_k \} $ of $\widetilde{\mathfrak{H}}$ as follows 
\begin{equation}
\begin{array}{rl}
\mathbf{R}( {\sigma }^0_{2k} )& = {\widetilde{\sigma }}_k, \, \, \mbox{for $k = 1, \ldots , K$}  \\
\rowspace \mathbf{R}({\sigma }^{+}_m + {\sigma }^{-}_m) & = {\widetilde{\sigma }}_{m} , \, \, 
\mbox{for $m \ge M$.} 
\end{array}
\label{eq-s3onenwdot}
\end{equation}
Recall that $2M = ${\tiny $\left\{ \begin{array}{rl} 
\hspace{-5pt} N+2, & \hspace{-5pt} \mbox{if $N$ is even} \\
\hspace{-5pt} N+1, & \hspace{-5pt} \mbox{if $N$ is odd} \end{array} \right. $} and 
$2K =${\tiny $\left\{ \begin{array}{cl} \hspace{-5pt} N, & 
\hspace{-5pt} \mbox{if $N$ is even} \\
\hspace{-5pt} N-1, &\hspace{-5pt} \mbox{if $N$ is odd.} \end{array} \right. $} 
To define shifting operators on ${\mathfrak{H}}^{\mathrm{even}}$ recall that equation 
(\ref{eq-threenew}) defines the shifting operators 
$\widetilde{\mathbf{b}}$ and $\widetilde{\mathbf{a}}$ on $\widetilde{\mathfrak{H}}$. We may lift the shifting operator $\widetilde{\mathbf{a}}$ to the shifting operator ${\mathbf{a}}^{\mathrm{even}}$ 
on ${\mathfrak{H}}^{\mathrm{even}}$ by setting 
\begin{subequations}
\begin{equation}
\begin{array}{rl}
\mathbf{R}{\mathbf{a}}^{\mathrm{even}}\, {\sigma }^0_{2k} & = 
\widetilde{\mathbf{a}}\, \mathbf{R}\, {\sigma }^0_{2k} \, \, \mbox{for $k \le K$} \\
\rowspace \mathbf{R}\, {\mathbf{a}}^{\mathrm{even}}({\sigma }^{+}_m+{\sigma }^{-}_m) & = 
\widetilde{\mathbf{a}} \, \mathbf{R} ({\sigma }^{+}_m+{\sigma }^{-}_m)\, \, \mbox{for $m \ge M$} 
\end{array}
\label{eq-s3onenwstara}
\end{equation} 
and lift the shifting operator $\widetilde{\mathbf{b}}$ to the shifting operator 
${\mathbf{b}}^{\mathrm{even}}$ on ${\mathfrak{H}}^{\mathrm{even}}$ by setting 
\begin{equation}
\begin{array}{rl}
\mathbf{R}{\mathbf{b}}^{\mathrm{even}}\, {\sigma }^0_{2k} & = 
\widetilde{\mathbf{b}}\, \mathbf{R}\, {\sigma }^0_{2k} \, \, \mbox{for $k \le K$} \\
\rowspace \mathbf{R}\, {\mathbf{b}}^{\mathrm{even}}({\sigma }^{+}_m+{\sigma }^{-}_m)  & = 
\widetilde{\mathbf{b}} \, \mathbf{R}({\sigma }^{+}_m+{\sigma }^{-}_m) \, \,   \mbox{for $m \ge M$.}
\end{array}
\label{eq-s3onenwstarb}
\end{equation}
\end{subequations}
If $k+1 \le K$, then 
\begin{displaymath}
\mathbf{R}{\mathbf{b}}^{\mathrm{even}}\, {\sigma }^0_{2k} = 
\widetilde{\mathbf{b}}\mathbf{R}\, {\sigma }^0_{2k} 
= \widetilde{\mathbf{b}}\, {\widetilde{\sigma }}_k = {\widetilde{\sigma }}_{k+1} = 
\mathbf{R}\, {\sigma }^0_{2k+2} .
\end{displaymath}
So ${\mathbf{b}}^{\mathrm{even}}{\sigma }^0_{2k} = {\sigma }^0_{2k+2}$, 
because $\mathbf{R}$ is injective. Since $n =2k$, it follows that ${\mathbf{b}}^{\mathrm{even}}$ 
raises the quantum number $n$ by $2$, provided that $n+1 < N$. Hence  
${\sigma }^0_{n+2} = \mathbf{b}{\sigma }^0_{n+1} = \mathbf{b}\mathbf{b} {\sigma }^0_n = 
{\mathbf{b}}^{\mathrm{even}}{\sigma }^0_n$.  A similar argument shows that ${\mathbf{b}}^{\mathrm{even}}$ raises the 
quantum number $m \ge M$ by $1$ and that 
${\mathbf{b}}^{\mathrm{even}}({\sigma }^{+}_m + {\sigma }^{-}_m) = 
{\mathbf{Q}}_{{\mathrm{e}}^{i\vartheta }}({\sigma }^{+}_m + {\sigma }^{-}_m)$. Analogous results can be obtained for 
the lowering operator ${\mathbf{a}}^{\mathrm{even}}$. In particular, if $k \le K$, then 
\begin{displaymath}
\mathbf{R}{\mathbf{a}}^{\mathrm{even}}{\sigma }^0_{2k} = 
\widetilde{\mathbf{a}}\, \mathbf{R}\, {\sigma }^0_{2k} = 
\widetilde{\mathbf{a}}\, {\widetilde{\sigma }}_k = {\widetilde{\sigma }}_{k-1} 
= \mathbf{R}\, {\sigma }^0_{2k-2}.
\end{displaymath}
This implies that ${\mathbf{a}}^{\mathrm{even}}{\sigma }^0_n = {\sigma }^0_{n-2} = 
{\mathbf{Q}}_{{\mathrm{e}}^{-2i\vartheta }}{\sigma }^0_n$ for $0 < n \le N$. Similarly, 
for $m > M$ we get ${\mathbf{a}}^{\mathrm{even}}({\sigma }^{+}_m + {\sigma }^{-}_m) = 
(m-1) ({\sigma }^{+}_m + {\sigma }^{-}_m)$.  

\section{Crossing the singularity}\label{section-7}
%%%%%%%%%%%%%%%%%%

The operators ${\mathbf{a}}^{\mathrm{even}}$ and ${\mathbf{b}}^{\mathrm{even}}$, defined in 
equation (\ref{eq-s3onenwstara}) and (\ref{eq-s3onenwstarb}), respectively allow for shifting quantum states which cross the singular level set $H^{-1}(2)$. 
In order to write this out explicitly, we need to consider the cases when $N$ is even or odd 
seperately. \medskip 

We look at the operators ${\mathbf{a}}^{\mathrm{even}}$ and ${\mathbf{b}}^{\mathrm{even}}$ 
when $N = 2K$, $M = \onehalf (N+2) = K+1$, and 
${\sigma }^0_{N}  = {\sigma }^0_{2K} \in {\mathfrak{H}}^{\mathrm{even}} \cap {\mathfrak{H}}_0$ 
together with  
\begin{align}
{\sigma }^{+}_M + {\sigma }^{-}_M & = {\sigma }^{+}_{K+1} +{\sigma }^{-}_{K+1} \in 
{\mathfrak{H}}^{\mathrm{even}} \cap ({\mathfrak{H}}_{+} \oplus {\mathfrak{H}}_{-}) .  \notag 
\end{align}
In this case equation (\ref{eq-s3onenwstarb}) yields 
\begin{align}
\mathbf{R}{\mathbf{b}}^{\mathrm{even}} {\sigma }^0_N & = 
\mathbf{R}{\mathbf{b}}^{\mathrm{even}}{\sigma }^0_{2K} = 
\widetilde{\mathbf{b}}\, \mathbf{R}\, {\sigma }^0_{2K} = 
\widetilde{\mathbf{b}}\, {\widetilde{\sigma }}_K \notag \\
& = {\widetilde{\sigma }}_{K+1} = {\widetilde{\sigma }}_M = 
\mathbf{R}({\sigma }^{+}_M + {\sigma }^{-}_M).
\end{align}
Since $\mathbf{R}: {\mathfrak{H}}^{\mathrm{even}} \rightarrow \widetilde{\mathfrak{H}}$ is injective, 
we get 
\begin{equation}
{\mathbf{b}}^{\mathrm{even}}  {\sigma }^0_N = {\sigma }^{+}_M + {\sigma }^{-}_M. 
\label{eq-ns7one}
\end{equation}
Let ${\mathbf{pr}}^{\pm}: {\mathfrak{H}}_{+} \oplus {\mathfrak{H}}_{-} \rightarrow {\mathfrak{H}}_{\pm}: 
{\sigma }^{+}_m + {\sigma }^{-}_m  \mapsto {\sigma }^{\pm}_m$. From (\ref{eq-ns7one}) we get 
\begin{equation}
{\mathbf{pr}}^{+}{\mathbf{b}}^{\mathrm{even}}  {\sigma }^0_N = {\sigma }^{+}_M \, \, \, 
\mathrm{and} \, \, \, 
{\mathbf{pr}}^{-}{\mathbf{b}}^{\mathrm{even}}  {\sigma }^0_N = {\sigma }^{-}_M,  
\label{eq-ns7two}
\end{equation}
which represents the transition given by the right pointing top and bottom slanted arrows in diagram (\ref{eq-news5two}) when $N$ is even. Similarly, 
\begin{align}
\mathbf{R}{\mathbf{a}}^{\mathrm{even}}({\sigma }^{+}_M + {\sigma }^{-}_M) & = 
\widetilde{\mathbf{a}}\, \mathbf{R}({\sigma }^{+}_M + {\sigma }^{-}_M) = 
\widetilde{\mathbf{a}}\, \mathbf{R}({\sigma }^{+}_{K+1} + {\sigma }^{-}_{K+1}) \notag \\
& = \widetilde{\mathbf{a}}\, {\widetilde{\sigma }}_{K+1} = {\widetilde{\sigma }}_K = 
\mathbf{R}\, {\sigma }^0_{2K} = \mathbf{R}\, {\sigma }^0_N,   \notag
\end{align}
which implies 
\begin{equation}
{\mathbf{a}}^{\mathrm{even}}({\sigma }^{+}_M + {\sigma }^{-}_M) = {\sigma }^0_N.
\label{eq-ns7three}
\end{equation}

Next we look at the operators ${\mathbf{a}}^{\mathrm{even}}$ and 
${\mathbf{b}}^{\mathrm{even}}$ when $N = 2K+1$ and $M = K+1$ and 
${\sigma }^0_N = {\sigma }^0_{2K+1} \in {\mathfrak{H}}^{\mathrm{odd}} \cap {\mathfrak{H}}_0$. 
Then by theorem 8 we have ${\sigma }^0_{N-1} = {\sigma }^0_{2K} \in 
{\mathfrak{H}}^{\mathrm{even}}\cap {\mathfrak{H}}_0$. Moreover, 
${\sigma }^{+}_M +{\sigma }^{-}_M = {\sigma }^{+}_{K+1}+{\sigma }^{-}_{K+1}  \in 
{\mathfrak{H}}^{\mathrm{even}} \cap ({\mathfrak{H}}_{+} \oplus {\mathfrak{H}}_{-})$. We can cross directly from ${\sigma }^0_{N-1} = {\sigma }^0_{2K}$ to 
${\sigma }^{+}_M + {\sigma }^{-}_M = {\sigma }^{+}_{K+1} + {\sigma }^{-}_{K+1}$ using the operator 
${\mathbf{b}}^{\mathrm{even}}$. In other words, 
\begin{equation}
{\mathbf{b}}^{\mathrm{even}}{\sigma }^0_{N-1} = {\sigma }^{+}_M + {\sigma }^{-}_M.
\label{eq-ns7four}
\end{equation}
Therefore, in order to cross from ${\sigma }^0_N$ to ${\sigma }^{+}_M + {\sigma }^{-}_M$, we first go to ${\sigma }^0_{N-1}$ and then to ${\sigma }^{+}_M + {\sigma }^{-}_M$. So 
${\mathbf{b}}^{\mathrm{even}} \mathbf{a} \, {\sigma }^0_N = {\sigma }^{+}_M + {\sigma }^{-}_M$. 
Hence for odd $N$, we have 
\begin{equation}
{\mathbf{pr}}^{+}{\mathbf{b}}^{\mathrm{even}} \mathbf{a} \, {\sigma }^0_N = {\sigma }^{+}_M 
\, \, \, \mathrm{and} \, \, \, 
{\mathbf{pr}}^{-}{\mathbf{b}}^{\mathrm{even}} \mathbf{a} \, {\sigma }^0_N = {\sigma }^{-}_M, 
\label{eq-ns7six}
\end{equation}
which represents the right pointing top and bottom slanted arrows in diagram 
(\ref{eq-news5two}) when $N$ is odd. Similarly, 
\begin{align}
\mathbf{R}{\mathbf{a}}^{\mathrm{even}}({\sigma }^{+}_M + {\sigma }^{-}_M) & = 
\widetilde{\mathbf{a}}\, \mathbf{R}({\sigma }^{+}_M + {\sigma }^{-}_M) = 
\widetilde{\mathbf{a}}\, \mathbf{R}({\sigma }^{+}_{K+1} + {\sigma }^{-}_{K+1})) \notag \\
& = \widetilde{\mathbf{a}}\, {\widetilde{\sigma }}_{K+1} = {\widetilde{\sigma }}_K = 
\mathbf{R}\, {\sigma }^0_{2K} = \mathbf{R}\, {\sigma }^0_{N-1}, \notag 
\end{align}
which implies 
\begin{equation}
{\mathbf{a}}^{\mathrm{even}}({\sigma }^{+}_M + {\sigma }^{-}_M) =  {\sigma }^0_{N-1} .
\label{eq-ns7seven}
\end{equation}
Let 
\begin{equation}
{\iota }_{\pm}: {\mathfrak{H}}_{\pm} \rightarrow {\mathfrak{H}}_{+} \oplus {\mathfrak{H}}_{-}: 
{\sigma }^{\pm}_m \mapsto \onehalf ({\sigma }^{\pm}_m + \mathbf{P}{\sigma }^{\pm}_m) = 
{\sigma }^{+}_m + {\sigma }^{-}_m. 
\label{eq-ns7threestar}
\end{equation}
Using the injection mapping ${\iota }_{\pm}$ and equations 
(\ref{eq-ns7three}) and (\ref{eq-ns7seven}) when $N = 2M -2$ we have 
\begin{subequations}
\begin{align}
{\widehat{\mathbf{a}}}^{\, \mathrm{even}}_{\pm} ({\sigma }^{\pm }_M) & = 
{\mathbf{a}}^{\mathrm{even}}{\iota }_{\pm}({\sigma }^{\pm }_M)
= {\mathbf{a}}^{\mathrm{even}}({\sigma }^{+}_M + {\sigma }^{-}_M)  = {\sigma }^0_N ;
\label{eq-ns7sevenstara}
\end{align}
while when $N = 2M -1$ we have 
\begin{align}
{\widehat{\mathbf{a}}}^{\, \mathrm{even}}_{\pm }({\sigma }^{\pm}_M) & = 
\mathbf{b}{\mathbf{a}}^{\mathrm{even}}{\iota }_{\pm}({\sigma }^{\pm}_M) 
 = \mathbf{b}{\mathbf{a}}^{\mathrm{even}}({\sigma }^{+}_M + {\sigma }^{-}_M)  
= \mathbf{b}{\sigma }^0_{N-1} = {\sigma }^0_N.
\label{eq-ns7sevenstarb}
\end{align}
\end{subequations} 
The operator ${\widehat{\mathbf{a}}}^{\, \mathrm{even}}_{+}$ represents  the transition given by 
the left pointing top slanting arrow in diagram (\ref{eq-news5two}); while the operator 
${\widehat{\mathbf{a}}}^{\, \mathrm{even}}_{-}$ represents  the transition given by 
the left pointing bottom slanting arrow in the diagram. 

\section{Appendix: construction of the lowering operator}
%%%%%%%%%%%%%%%%%%%%%%%

In this appendix we construct the lowering operator $\mathbf{a}$ for the quantized 
${\Z }_2$-reduced system on $T^{\times} S^1$. \medskip 

On phase space $T^{\ast }{\widetilde{S}}^1 = \R \times {\widetilde{S}}^1$, where 
${\widetilde{S}}^1 = \R/(2\pi \, \Z)$ with coordinates $(p, \vartheta )$ and symplectic form 
$\widetilde{\omega} = \dee p \wedge \dee \vartheta $ consider the trivial (right) principal bundle 
\begin{displaymath}
{\pi }^{\times }: L^{\times } = {\C }^{\times } \times T^{\ast }{\widetilde{S}}^1 \rightarrow 
T^{\ast }{\widetilde{S}}^1: \big( b, (p, \vartheta ) \big) \mapsto (p, \vartheta )
\end{displaymath}
with connection $1$-form $\beta = \frac{1}{2\pi } \frac{\dee b}{b} - \frac{1}{h} p\, \dee \vartheta $. The curvature $\dee \beta $ of $\beta $ is $-\frac{1}{h} \widetilde{\omega }$, which we suppose has integer de Rham cohomology. \medskip 

The action integral of Bohr-Sommerfeld quantization is $I = \int^{2\pi }_0 p\, \dee \vartheta = 2\pi \, p$. The variable conjugate to $I$ is $\theta = \vartheta /2\pi $, since 
$\widetilde{\omega} = \dee I \wedge \dee \theta $. Let $(I, \theta)$ be coordinates on 
$T^{\ast }S^1 = \R \times S^1$ with $S^1 = \R /\Z $ and let 
$\omega  = \dee I \wedge \dee \theta $ be the symplectic form on $T^{\ast }S^1$. 
The vector field $X = -\frac{\partial }{\partial I} = 
-\frac{1}{2\pi } \frac{\partial }{\partial p}$ on $T^{\ast }S^1$ is locally Hamiltonian, since 
\begin{displaymath}
L_X\omega = \dee (X \lefthook \omega ) + X \lefthook \dee \omega = \dee (-\dee \theta ) =0, 
\end{displaymath}
and has \emph{local} Hamiltonian $\theta $, since $-\dee \theta = X \lefthook \omega $ 
\emph{locally}. The flow of $X$ is 
\begin{equation}
{\mathrm{e}}^{tX}: T^{\ast }S^1 \rightarrow T^{\ast }S^1: (I, \theta ) \mapsto 
(I-t, \theta ) = \big( 2\pi (p - \ttfrac{1}{2\pi }t) , \vartheta /2\pi \big) . 
\label{eq-exone}
\end{equation}
Note that the diffeomorphism ${\mathrm{e}}^{hX}$ sends the Bohr-Sommerfeld 
torus $T_n$, defined by $\{ I = nh \} $, onto the Bohr-Sommerfeld torus 
$T_{n-1}$, defined by $\{ T = (n-1)h \}$. \medskip 

In what follows we find a quantomorphism ${\Phi}_h$ of $(L, \lambda )$, which covers 
${\mathrm{e}}^{hX}$. In other words, ${\Phi }_h$ is a diffeomorphism of $L$ into itself 
such that ${\Phi }^{\ast }_h\lambda = \lambda $ and $\pi \comp {\Phi }_h = 
{\mathrm{e}}^{hX} \comp \pi $. Here 
\begin{displaymath}
\pi : L = \C \times T^{\ast }S^1 \rightarrow T^{\ast }S^1:\big( z, (I, \theta ) \big) \mapsto 
(I, \theta ) 
\end{displaymath}
is the line bundle, associated to the ${\C }^{\times }$ principal bundle 
${\pi }^{\times }: L^{\times } \rightarrow T^{\ast}{\widetilde{S}}^1$, with connection $1$-form $\lambda = \frac{1}{2\pi i} \dee z - \frac{1}{h} I \dee \theta $. \medskip 

The vector field $X$ is \emph{integral}, that is, 1) there is an \emph{good covering} 
$\mathcal{U} = {\{ U_i \} }_{i \in I}$ of $T^{\ast }S^1$ by open sets $U_i$, $i \in I$, where 
every finite intersection of elements of $\mathcal{U}$ is either empty or contractible; 
2) for every $U_i$, $U_j \in \mathcal{U}$ such that $U_i \cap U_j \neq \varnothing $ we 
have $\theta |_{U_i} - \theta |_{U_j}$ is an integer on $U_i \cap U_j$. The local 
Hamiltonian functions $\theta |_{U_i}$ for $i \in I$ piece together to give a 
smooth mapping $[\theta ]:T^{\ast }S^1 \rightarrow S^1 = \R /\Z $, which is the 
``coordinate'' $\theta $, that is, $[\theta ] = \vartheta /2\pi $. \medskip 

Consider the vector field $Z$ on $L^{\times }$, whose flow is 
\begin{equation}
{\mathrm{e}}^{tZ}: L^{\times } \rightarrow L^{\times }: \big( b, (I, \theta ) \big) 
\mapsto \big( b\, {\mathrm{e}}^{-2\pi i \, t \theta /h}, (I-t, \theta ) \big) . 
\label{eq-extwo}
\end{equation}
The flow of $Z$ preserves the connection $1$-form $\beta $, since 
for every $(b, I ,\theta ) ) \in L^{\times }$ we have 
\begin{align}
\big( ({\mathrm{e}}^{tZ})^{\ast }\beta \big) (b,I,\theta ) & = 
\ttfrac{1}{2\pi i} \dee \ln [ b \, {\mathrm{e}}^{-2\pi i\, t \theta /h}]
 - \ttfrac{1}{h} (I-t) \dee \theta \notag \\
& = \big( \ttfrac{1}{2\pi i} \frac{\dee b}{b} - \ttfrac{1}{h} I \dee \theta \big) 
-\ttfrac{1}{h}t \, \dee \theta + \ttfrac{1}{h}t \, \dee \theta = \beta (b, I, \theta ). \notag 
\end{align} 
We have 
\begin{equation}%
{\mathrm{e}}^{tZ} = {\mathrm{e}}^{-tY_{\theta /h}} \comp {\mathrm{e}}^{t \, \mathrm{lift}X}, 
\label{eq-exthree}
\end{equation}
where $Y_{\theta /h}(b, I, \theta ) = \smalldbydt (b \, {\mathrm{e}}^{2\pi i \, t\theta /h}, I, \theta )$ 
is a vector field on $(L^{\times }, \beta )$, whose flow is 
${\mathrm{e}}^{tY_{\theta /h}}(b,I,\theta ) = (b\, {\mathrm{e}}^{2\pi i\, t\theta /h}, I, \theta )$, and 
$\mathrm{lift}X$ is a vector field on $(L^{\times }, \beta )$, which is the horizontal lift 
of the vector field $X$, that is, $\mathrm{lift}X(b,I, \theta ) \in \ker \beta (b,I, \theta )$ for 
every $(b,I,\theta ) \in L^{\times }$. The vector fields $\mathrm{lift}X$ and $X$ are 
${\pi }^{\times }$-related, that is, 
$T_{(b,I,\theta )}{\pi }^{\times } \big( \mathrm{lift}X (b,I,\theta ) \big) 
= X\big( {\pi }^{\times }(b,I,\theta ) \big) $ for every $(b,I,\theta ) \in L^{\times }$. The 
flow of $\mathrm{lift} X$ is 
\begin{displaymath}
{\mathrm{e}}^{t\, \mathrm{lift}X}(b,I, \theta ) = \big( b, {\mathrm{e}}^{t X}(I,\theta ) \big) = 
(b, I-t , \theta ).
\end{displaymath}
Note that the flows of the vector fields $Y_{\theta /h}$ and $\mathrm{lift} X$ commute. \medskip 

We now look at the universal covering space $(T^{\ast }\R , \widetilde{\omega })$ of 
$(T^{\ast }S^1, \omega )$ with coordinates $(p,q)$ and symplectic form 
$\widetilde{\omega } = \dee p \wedge \dee q$. The universal covering map is given by 
\begin{displaymath}
\kappa : T^{\ast }\R \rightarrow T^{\ast }S^1: (p,q) \mapsto (\ttfrac{1}{2\pi }I, \theta ) = 
(p, q \bmod \Z ) , 
\end{displaymath}
since $\R \rightarrow S^1 = \R /\Z : q \mapsto q \bmod \Z $ is the smooth universal 
covering map of $S^1$. Pull the local Hamiltonian vector field $X$ on $T^{\ast }S^1$ 
back by the covering map $\kappa $ to a vector field $\widetilde{X}$ on $T^{\ast }\R $, 
which is $\kappa $-related to $X$, that is, $T_{(p,q)}\kappa \big( \widetilde{X}(p,q) \big) = 
X\big( \kappa (p,q) \big)$ for every $(p,q) \in T^{\ast }\R$. 
The integral vector field $\widetilde{X}$ 
is the Hamiltonian vector field $X_{q+c} = - \frac{\partial}{\partial p}$ associated to 
the Hamiltonian function $q:T^{\ast }\R \rightarrow \R: (p,q) \mapsto q$. The constant 
$c$ can be choosen to be $0$, because the smooth mappings $[q]: T^{\ast }\R \rightarrow S^1 
= \R /\Z : (p, q) \mapsto q \bmod \Z $ and ${\kappa }^{\ast }[\theta ]: T^{\ast }\R 
\rightarrow S^1: (q,p) \mapsto [\theta ]\big( \kappa (p,q) \big) $ are equal, namely, 
\begin{equation}
[q] = {\kappa }^{\ast }[\theta ]. 
\label{eq-exthreestar}
\end{equation}
A calculation shows that 
\begin{subequations}
\begin{equation}
\kappa \comp {\mathrm{e}}^{tX_q} = {\mathrm{e}}^{tX} \comp \kappa . 
\label{eq-exthreea}
\end{equation}
\end{subequations}
The ${\C }^{\times }$ bundle ${\pi }^{\times} : L^{\times } \rightarrow T^{\ast }S^1$ pulls 
back under the covering map $\kappa $ to the $C^{\times }$ bundle 
${\widetilde{\pi }}^{\times }: {\widetilde{L}}^{\times } \rightarrow T^{\ast }\R : 
(b, p,q ) \mapsto (p,q)$ with connection $1$-form $\widetilde{\beta } = {\kappa }^{\ast }\beta = 
\ttfrac{1}{2\pi i} \frac{\dee b}{b} - \ttfrac{1}{h} p \, \dee q$. The flow ${\mathrm{e}}^{tZ}$ on 
$L^{\times }$ lifts to the $1$-parameter group of diffeomorphisms 
\begin{equation}
{\mathrm{e}}^{t\widetilde{Z}}: {\widetilde{L}}^{\times } \rightarrow {\widetilde{L}}^{\times }: 
(b, p, q) \mapsto \big( b\, {\mathrm{e}}^{-2\pi i \, tq/h}, p-t ,q \big) , 
\label{eq-exfour}
\end{equation}
which preserves $\widetilde{\beta }$. We have 
${\mathrm{e}}^{t\widetilde{Z}} = {\mathrm{e}}^{-t{\widetilde{Y}}_{q/h}} \comp 
{\mathrm{e}}^{t\, \mathrm{lift}X_q}$, 
where ${\widetilde{Y}}_{q/h}$ is the vector field on ${\widetilde{L}}^{\times }$, whose 
value at $(b,p,q)$ is $\smalldbydt (b\, {\mathrm{e}}^{-2\pi i \, tq/h}, p,q)$. Its flow is given by 
${\mathrm{e}}^{tY_{q/h}}(b,p,q) = 
(b\, {\mathrm{e}}^{2\pi i \, tq/h}, p,q)$. Also $\mathrm{lift}X_q$ is a vector field on 
${\widetilde{L}}^{\times }$, which is the horizontal lift of $X_q$ using the connection 
$1$-form $\widetilde{\beta }$ on ${\widetilde{L}}^{\times }$. The flow of 
$\mathrm{lift}X_q$ is ${\mathrm{e}}^{t\, \mathrm{lift}X_q}(b,p,q) = (b, p-t ,q)$. Note that 
the flows ${\mathrm{e}}^{tY_{q/h}}$ and ${\mathrm{e}}^{t\, \mathrm{lift}X_q}$ commute. 
Since ${\widetilde{L}}^{\times } = {\kappa }^{\ast }L^{\times }$, there is a smooth 
mapping 
\begin{displaymath}
{\kappa }^{\times }: {\widetilde{L}}^{\times } \rightarrow L^{\times }: (b,p,q) 
\mapsto (b, p, q \bmod \Z ) = (b, \ttfrac{1}{2\pi }I, \theta ) , 
\end{displaymath}
which covers $\kappa $, that is, ${\pi }^{\times } \comp {\kappa }^{\times } = 
\kappa \comp {\widetilde{\pi }}^{\times }$. The flows ${\mathrm{e}}^{t \, \mathrm{lift}X_q}$  and 
${\mathrm{e}}^{t\, \mathrm{lift}X}$ are ${\kappa }^{\times }$-related, that is 
\addtocounter{equation}{-2}
\begin{subequations}
\addtocounter{equation}{1}
\begin{equation}
{\kappa }^{\times } \comp {\mathrm{e}}^{t \, \mathrm{lift}X_q} = 
{\mathrm{e}}^{t \, \mathrm{lift}X} \comp {\kappa }^{\times }. 
\label{eq-exthreeb}
\end{equation}
\end{subequations}
 
Let 
\begin{displaymath}
{\sigma }^{\times }: T^{\ast }S^1 \rightarrow L^{\times}:
(I, \theta ) \mapsto \big( b(I,\theta ), I, \theta \big) 
\end{displaymath}
be a smooth section of the bundle ${\pi }^{\times }: L^{\times } \rightarrow T^{\ast }S^1$, 
where $(I,\theta ) \mapsto b(I,\theta )$ is a smooth nowhere vanishing complex valued 
function on $T^{\ast }S^1$. Then 
\begin{align}
({\mathrm{e}}^{tZ})_{\ast }{\sigma }^{\times } & = 
{\mathrm{e}}^{tZ} \comp {\sigma }^{\times } \comp {\mathrm{e}}^{-tX} 
= {\mathrm{e}}^{-tY_{\theta /h}} \comp  
\big( {\mathrm{e}}^{t\, \mathrm{lift}X} \comp {\sigma }^{\times } \comp {\mathrm{e}}^{-tX} \big) 
\notag \\
& = \big( {\mathrm{e}}^{-tY_{\theta /h}} \big)_{\ast }  \comp 
\big( {\mathrm{e}}^{t\, \mathrm{lift}X} \big)_{\ast }{\sigma }^{\times } 
\label{eq-exfive}
\end{align}  
is a smooth section of the bundle ${\pi }^{\times }: L^{\times } \rightarrow T^{\ast }S^1$. 
Let 
\begin{displaymath}
{\widetilde{\sigma }}^{\times }: T^{\ast }\R \rightarrow {\widetilde{L}}^{\times }:
(p,q) \mapsto \big( \, \widetilde{b}(p,q), p,q\big) 
\end{displaymath}
be a smooth section of the bundle ${\widetilde{\pi }}^{\times }: {\widetilde{L}}^{\times } 
\rightarrow T^{\ast }\R$, which is the pull back by the mapping ${\kappa }^{\times }$ of 
the smooth section ${\sigma }^{\times }$ of the bundle ${\pi }^{\times }: L^{\times } 
\rightarrow T^{\ast }S^1$. For every $(p,q) \in T^{\ast }\R$ we  have 
${\widetilde{\sigma }}^{\times }(p,q) = \big( b(\kappa (p,q)), (p,q) \big)$. So 
\addtocounter{equation}{-2}
\begin{subequations}
\addtocounter{equation}{2}
\begin{equation}
{\kappa }^{\times } \comp {\widetilde{\sigma }}^{\times } = {\sigma }^{\times} \comp \kappa , 
\label{eq-exthreec}
\end{equation}
\end{subequations}
which characterizes ${\widetilde{\sigma }}^{\times }$. We now show 
that
\begin{equation}
({\mathrm{e}}^{h\widetilde{Z}})_{\ast }{\widetilde{\sigma }}^{\times } 
=({\kappa }^{\times })^{\ast } \big( {\mathrm{e}}^{-2 \pi i \, [\theta]} \mbox{\tiny $\bullet$} 
({\mathrm{e}}^{h\, \mathrm{lift}X })_{\ast }{\sigma }^{\times } \big) , 
\label{eq-exsix}
\end{equation}
We explain what the symbol \mbox{\tiny $\bullet$} in the above formula means. For a smooth section ${\sigma }^{\times}:T^{\ast }S^1 \rightarrow L^{\times }:
(I, \theta ) \mapsto \big( b(I,\theta ), I, \theta \big)$ of the bundle 
${\pi }^{\times }:L^{\times } \rightarrow T^{\ast }S^1$ and a smooth complex 
valued function $f$ on $T^{\ast }S^1$ we define $f \mbox{\tiny $\bullet$}\, {\sigma }^{\times }: 
T^{\ast }S^1 \rightarrow L^{\times }$ to be the smooth section $(I, \theta ) 
\mapsto \big( b(I,\theta )f(I,\theta ), I, \theta \big)$. \medskip  

Analgous to (\ref{eq-exthree}) we have 
\begin{equation}
\big( {\mathrm{e}}^{t\widetilde{Z}} \big)_{\ast }{\widetilde{\sigma }}^{\times }  = 
\big( {\mathrm{e}}^{-t{\widetilde{Y}}_{\theta /h}} \big)_{\ast } \comp 
\big( {\mathrm{e}}^{t \, \mathrm{lift}X_q} \big)_{\ast } {\widetilde{\sigma }}^{\times }, 
\label{eq-exseven}
\end{equation}
where ${\widetilde{\sigma }}^{\times } = ({\kappa }^{\times })^{\ast }{\sigma }^{\times }$. \medskip 

\noindent \textbf{Verification of (70)} To 
verify equation (\ref{eq-exsix}) we note that 
\addtocounter{equation}{-1}
\begin{subequations}
\begin{equation}
({\kappa }^{\times })^{\ast } \big( ({\mathrm{e}}^{h\, \mathrm{lift}X})_{\ast }{\sigma }^{\times } \big)
= ({\mathrm{e}}^{h\, \mathrm{lift}X_q})_{\ast } \big( 
({\kappa }^{\times })^{\ast }{\sigma }^{\times } \big)  
\label{eq-exsevena}
\end{equation}
follows by applying equations (\ref{eq-exthreea}), (\ref{eq-exthreeb}), and 
(\ref{eq-exthreec}). Next we show that 
\begin{equation}
({\mathrm{e}}^{h{\widetilde{Y}}_{q/h}})_{\ast } 
\big( ({\kappa }^{\times })^{\ast }{\sigma }^{\times } \big) = 
({\kappa }^{\times })^{\ast} (
{\mathrm{e}}^{-2\pi i \, [\theta ]} \mbox{\tiny $\bullet $}\, {\sigma }^{\times }) .
\label{eq-exsevenb}
\end{equation}
\end{subequations}
Using (\ref{eq-exthreestar}), for every $(p,q) \in T^{\ast }\R $ we get 
\begin{align} 
{\mathrm{e}}^{2\pi i \, [\theta ](\kappa (p,q)) } & 
= {\mathrm{e}}^{2\pi i \, [q]} = 
{\mathrm{e}}^{2\pi i (q + n)}, \, \, \mbox{for every $n \in \Z $} \notag \\
& = {\mathrm{e}}^{2\pi i \, q} = {\mathrm{e}}^{h(2\pi i \, q/h)}. \notag 
\end{align}
So 
\begin{align}
({\mathrm{e}}^{h{\widetilde{Y}}_{q/h}})_{\ast } 
\big( ({\kappa }^{\times })^{\ast }{\sigma }^{\times } \big) (p,q) & = 
({\kappa }^{\times})^{\ast }\big( ({\mathrm{e}}^{hY_{q/h}})_{\ast }{\sigma }^{\times } \big) (p,q)   
= \big( ({\mathrm{e}}^{h Y_{[\theta ]/h}})_{\ast }{\sigma }^{\times } \big) \big( \kappa (p,q) \big) \notag \\
&\hspace{-1.4in} = 
{\sigma }^{\times }\big( \kappa (p,q) \big) {\mathrm{e}}^{h(2\pi i \, q/h)} \
= {\sigma }^{\times }\big( \kappa (p,q) \big) {\mathrm{e}}^{2\pi i \, [\theta ](\kappa (p,q))} 
 = \big( {\mathrm{e}}^{2\pi i\, [\theta ]} \mbox{\tiny $\bullet$}\, {\sigma }^{\times }\big) 
\big( \kappa (p,q) \big)  \notag \\
& \hspace{-1.4in}
 = ({\kappa }^{\times })^{\ast }({\mathrm{e}}^{2\pi i\, [\theta ]} \mbox{\tiny $\bullet$}\, 
{\sigma }^{\times } \big) (p,q) .   \notag
\end{align}
Equation (\ref{eq-exseven}) follows from equations (\ref{eq-exsevena}) and 
(\ref{eq-exsevenb}). \quad \mbox{\tiny $\blacksquare $} \medskip 

The locally Hamiltonian vector field $X$ on $T^{\ast }S^1$ with flow 
${\mathrm{e}}^{tX}:T^{\ast }S^1 \rightarrow T^{\ast }S^1: (I, \theta ) \mapsto (I-t,\theta )$ 
lifts to a vector field $\widehat{X}$ on $L$, whose flow is 
${\mathrm{e}}^{t\widehat{X}}: L \rightarrow L: (z, I, \theta ) \mapsto (z, I-t, \theta )$. The map 
\begin{equation}
\mu : L = \C \times T^{\ast }S^1 \rightarrow L^{\times } = {\C }^{\times } \times T^{\ast }S^1:
(z, I, \theta ) \mapsto ({\mathrm{e}}^{z}, I, \theta ) = (b, I, \theta )
\label{eq-exeight}
\end{equation}
is smooth and ${\mu }^{\ast }\beta = \lambda $. Moreover, we have 
$\mu \comp {\mathrm{e}}^{t\widehat{X}} = {\mathrm{e}}^{t\, \mathrm{lift}X} \comp \mu $. 
Instead of ${\mathrm{e}}^{t\widehat{X}}$ we will write ${\widehat{\mathrm{e}}}^{\, tX}$. 
The mapping $\mu $ intertwines the (right) action of ${\C }^{\times }$ on $L^{\times }$ 
with the action of ${\C }^{\times }$ on $L$, namely, $\mu (b'z, I, \theta ) = 
\mu (b,I,\theta )(b')^{-1}$ for every $b$, $b' \in {\C }^{\times }$ and every $(I,\theta ) \in 
T^{\ast }S^1$. From these remarks it follows that the operator ${\Phi }^{\times }_h = 
{\mathrm{e}}^{-2\pi i\, [\theta ]} \mbox{\tiny $\bullet $}\, 
({\widehat{\mathrm{e}}}^{\, h \, \mathrm{lift}X})_{\ast }$ on smooth sections of the line 
bundle ${\pi }^{\times }: L^{\times } \rightarrow T^{\ast }S^1$ becomes the operator 
${\Phi }_h = {\mathrm{e}}^{2\pi i \, [\theta ]} 
({\widehat{\mathrm{e}}}^{\, h\, \mathrm{lift}X})_{\ast }$ 
on smooth sections of the line bundle $\pi : L \rightarrow T^{\ast }S^1$. Clearly, 
${\Phi }_h$ covers ${\mathrm{e}}^{hX}$, that is, $\pi \comp {\Phi }_h = {\mathrm{e}}^{hX} 
\comp \pi $, and preserves the connection $1$-form $\lambda $. So ${\Phi }_h$ is the desired lowering operator $\mathbf{a}$. We note that ${\Phi }_{-h}$ is the raising operator $\mathbf{b}$.


\begin{thebibliography}{99}

\bibitem{bohr} Niels Bohr, On the constitution of atoms and molecules (Part I) , \textit{Philosophical Magazine}, \textbf{26} (1913) 1-25.

\bibitem{cushman-bates15} Richard H. Cushman and Larry M. Bates, ``Global aspects of classical integrable systems'', second edition, Birkh\"{a}user, Basel, 2015. 

\bibitem{cushman-sniatycki13} Richard Cushman and J\k{e}drzej \'{S}niatycki, 
Bohr-Sommerfeld-Heisenberg theory in geometric quantization, \textit{J.
Fixed Point Theory Appl.} \textbf{13} (2013) 3--24.

\bibitem{cushman-sniatycki12} Richard Cushman and J\k{e}drzej \'{S}niatycki,
Bohr-Sommerfeld Heisenberg quantization of the $2$-dimensional harmonic
oscillator, \texttt{arXiv:math.SG.1207.1477v2}.

\bibitem{cushman-sniatycki13} Richard Cushman and J\k{e}drzej \'{S}niatycki,
Bohr-Sommerfeld-Heisenberg theory in geometric quantization, \textit{J.
Fixed Point Theory Appl.} \textbf{13} (2013) 3--24.

\bibitem{cushman-sniatycki18} Richard Cushman and J\k{e}drzej \'{S}niatycki, 
Shifting operators in geometric quantization, \texttt{arXiv:math.SG.1808.04002}.

\bibitem{dirac25} Paul Adrien Maurice Dirac, The fundamental equations of quantum
mechanics, \textit{Proc. Roy. Soc. London}, \textbf{A 109} (1925) 642--653.

\bibitem{dirac30} Paul Adrien Maurice Dirac, 
``The principles of quantum mechanics'', Clarendon Press, Oxford, UK, 1930

\bibitem{dullin} Holger Dullin, Semi-global symplectic invariants of the
spherical \linebreak pendulum, \textit{J. Differential Equations} \textbf{254} (2013)
2942--2963.

\bibitem{heisenberg25} Werner Heisenberg, \"{U}ber die quantentheoretische Umdeutung 
kinematischer und mechanischer Beziehungen, \textit{Z. Phys.} \textbf{33} (1925) 879--893.

\bibitem{sniatycki80} J\k{e}drzej \'{S}niatycki, \textit{Geometric
Quantization and Quantum Mechanics}, Applied Mathematical Series \textbf{30}, 
Springer Verlag, New York, 1980.

\bibitem{sommerfeld} Arnold Sommerfeld, Zur Theorie der Balmerschen Serie, 
\textit{Sitzungberichte der Bayerischen Akademie der Wissenschaften (M\"{u}%
nchen), mathematisch-physikalische Klasse}, (1915) 425-458.

\end{thebibliography}
\end{document}